\magnification 1200

\def\ver{{\it Version 1.00}}
\def\vdate{ {\it 4 March  2005}}

   %%SPACINGS  
\def\bigb{\bigskip\noindent}         
\def\med{\medskip\noindent}
\def\sms{\smallskip}
\def\ms{\medskip}
\def\bs{\bigskip}
\def\small{\smallskip\noindent}
\def\sm{\smallskip\noindent}
\def\nl{\hfil\break}
\def\noi{\noindent}

  %%TODAY'S DATE
\def\today{\noindent\number\day
\space\ifcase\month\or
  January\or February\or March\or April\or May\or June\or
  July\or August\or September\or October\or November\or December\fi
  \space\number\year}

  %%FORMATTING
\def\cl {\centerline}

\def\section#1#2{{\bigbreak\bigskip \centerline{\bf #1. #2}\bigskip}}
\def \nr { \smallskip \noindent }

  %%MATHS
\input amssym.def
\input amssym.tex

  %%Open
\def\bP {{\Bbb P}} \def\bE {{\Bbb E}} 
\def\bR {{\Bbb R}} \def\bN {{\Bbb N}} \def\bZ {{\Bbb Z}}
\def\bB {{\Bbb B}}

  %%Script
  \def\sC {{\cal C}}
 \def\sE {{\cal E}} \def\sF {{\cal F}}
\def\sG {{\cal G}}  
  \def\sL {{\cal L}}

  %% Overline

%\def\ob{\overline B}

  %%Greek
\def\al {\alpha}
\def\lam {\lambda}  \def\Gam{\Gamma}          
   \def\gam{\gamma}                

\def\ep{\varepsilon}  \def\eps{\varepsilon}
\def\om{\omega }
\def\th{\theta} \def\Th{\Theta}

  %%Maths

\def\pd {\partial}       
\def\q{\quad} \def\qq{\qquad}
\def\dint{\int\kern-.6em\int}

  %%Mathops

%\def\liminf{\mathop{\underline{\rm lim}}}
%\def\limsup{\mathop{{\rm lim\ sup}}}

  %%Fractions
\def \fract#1#2{{\textstyle {#1\over #2}}}
\def \frac#1#2{{ {#1\over #2}}}

\def \half {{\textstyle {1 \over 2}}}

\def \qed {\hfill$\square$\par}

%% Special

\def\=d{{\,\buildrel (d) \over =\,}}
\def\a.s.{{\buildrel a.s. \over \longrightarrow}}

\def\wt{\widetilde}

\def\proof{\sm {\it Proof. }}
\def\tf{\wt f}

\def\tbP {{\Bbb P}} \def\tbE {{\Bbb E}}
\def\Gam{\Gamma}  \def\gam{\gamma}

\def\eqd{{\,\buildrel (d) \over =\,}}
\def\wX{\wt X}
\def\wY{\wt Y}
\def\Bin{{\tt Bin}}
\def\Ber{{\tt Ber}}
\font\tf=cmbx12 scaled\magstephalf

\cl {\tf Random walk on the incipient infinite cluster on trees }

\vskip 0.3 truein

\centerline {Martin T. Barlow\footnote{$^1$}{Research partially
 supported by a grant from NSERC (Canada).},
\qquad  
Takashi Kumagai\footnote{$^2$}{Research partially supported by the 
Grant-in-Aid for Scientific Research for Young Scientists (B) 16740052.}}

\bigskip

{\narrower 
\med \noi{\bf Abstract.}
Let $\sG$ be the incipient infinite cluster
(IIC) for percolation on a homogeneous tree of degree $n_0+1$.
We obtain estimates for the transition density of the
the continuous time simple random walk 
$Y$ on $\sG$; the process satisfies anomalous diffusion 
and has spectral dimension $\fract43$.
}

\bigb
\noi{\bf 2000 MSC.} 
Primary 60K37; Secondary 60J80, 60J35.

\noi{\bf Keywords.} Percolation, incipient infinite cluster, 
random walk, branching process, heat kernel.

\bigskip
\med
{\bf 1. Introduction }
\ms
We recall the bond percolation model on the lattice $\bZ^d$: 
each bond is open with probability $p\in (0,1)$, independently
of all the others. Let $\sC(x)$ be the open cluster containing
$x$; then if $\th(p)=P_p(|\sC(x)|=+\infty)$ it is well
known (see [Gm]) that there exists $p_c=p_c(d)$ such that
$\th(p)=0$ if $p<p_c$ and $\th(p)>0$ if $p>p_c$. 

If $d=2$ or $d\ge 19$ (or $d > 6$ for `spread out' models) it is known
(see [Gm], [HS]) that $\th(p_c)=0$, and it is conjectured that this
holds for all $d\ge 2$. At the critical probability $p=p_c$ it is 
believed that in any box of side $n$ there exist with high probability
open clusters of diameter of order $n$ -- see [BCKS].  For large $n$ the
local properties of these large finite clusters can, in certain
circumstances, be captured by regarding them as subsets of an infinite
cluster $\wt \sC$, called the `incipient infinite cluster' (IIC).

This was constructed when $d=2$ in [Ke1], by taking the limit 
as $N \to \infty$ of the cluster $\sC(0)$ conditioned to intersect
the boundary of a box of side $N$ with center at the origin.
See [Ja1], [Ja2] for other constructions of the IIC in two 
dimensions.
For large $d$ a construction of the IIC in $\bZ^d$ is 
given in [HJ], using the lace expansion. 
It is believed that the results there will hold for any $d>6$. [HJ] also 
gives the existence and some properties of the IIC
for all $d>6$ for `spread-out' models: these include the case when
there is a bond between $x$ and $y$ with probability $p L^{-d}$
whenever $y$ is in a cube side $L$ with center $x$, and the parameter
$L$ is large enough.
Rather more is known about the IIC for oriented percolation
on $\bZ_+ \times \bZ^d$ 
(see [HHS], [HS]),
but in this discussion, which mainly concerns what is conjectured 
rather than what is known, we specialize to the case of $\bZ^d$.
We write $\wt\sC_d$ for the IIC in $\bZ^d$.
It is believed that the global properties of  $\wt\sC_d$ are the same for
all $d>d_c$, both for nearest neighbour and spread-out models.  In [HJ]
it is proved for `spread-out' models that $\wt\sC_d$ has one end --
that is that any two paths from 0 to infinity intersect infinitely often.

For large $d$, it is believed that the geometry of $\wt\sC_d$ is also
similar to that of the IIC when `$d=\infty$' --
that is to the IIC on a regular tree; 
this is supported by the results in [HHS] and [HJ]. For trees the
construction of the IIC is much easier than for lattices, and there is
a close connection between the IIC and a critical
Bienaym\'e-Galton-Watson branching processes conditioned on
non-extinction.  
In [Ke2] Kesten gave the construction of the IIC 
$\sG$ for critical branching processes. This is an infinite subtree,
which contains only one path from the root to infinity.
This tree is quite sparse, and has polynomial volume growth:
in the case when the offspring distribution has finite variance,
a ball $B(x,r)$ in $\sG$ has roughly $r^2$ points. (This is when
distance in $\sG$ is measured using the natural graph distance).

Let $Y=(Y_t, t\ge 0)$ be the simple random walk on $\wt\sC_d$,
and $q_t(x,y)$ be its transition density (see Section 3 for a
precise definition).
Define the {\sl spectral dimension} of  $\wt\sC_d$ by
$$ d_s(\wt\sC_d) =-2 \lim_{t \to \infty} \frac{\log q_t(x,x)}{\log t}, 
\eqno (1.1) $$
(if this limit exists). Alexander and Orbach [AO] conjectured that, 
for any $d\ge 2$, 
$d_s(\wt\sC_d)=4/3$. While it is now thought that this is unlikely
to be true for small $d$, the results on the geometry of  $\wt\sC_d$
in [HHS] and [HJ] are consistent with this holding for large $d$.
(Or for any $d$ above the critical dimension for spread-out models).

Random walks on supercritical clusters in $\bZ^d$
are studied in [B2] (transition density estimates) and [SS] (invariance
principle for $d\ge 4$). In these cases the large scale behaviour
of the random walk approximates that of the random walk on $\bZ^d$, 
and the unique infinite cluster has spectral dimension $d$.
 
In what follows, we will specialize to the case of critical 
percolation on a regular rooted tree with degree $n_0+1$, which
we denote $\bB$. We write $0$ for the root of $\bB$.
We keep $n_0$ fixed, but (in view of possible future applications) 
wish to obtain estimates which do not depend on $n_0$. 
For bond percolation with probability $p$ on $\bB$, it is easy to 
see that if $X_n$ is the 
number of vertices at level $n$ in $\sC(0)$, then $X=(X_n)$ is 
a branching process with $\Bin(n_0,p)$ offspring distribution.
Thus $p_c=1/n_0$. For the construction of the IIC see [Ke2]: 
we obtain a subtree $\sG\subset \bB$ with law $\bP$, on a 
probability space $(\Omega_1, \sF, \bP)$.
Write $\bB_N$ for the $N$-th level of $\bB$, and 
$\bB_{\le N}$ for the union of the first $N$ levels
of $\bB$. 
Then the law of $\sG$ is characterized by the
fact that the law of $\sG \cap \bB_{\le N}$ under $\bP$ is the same
as that of $\sC(0)$ under $P_{p_c}$, conditioned on $\sC(0)$
reaching level $N$. 

Motivated by [AO], in [Ke2] Kesten studied the simple random walk 
on $\sG(\om)$, and also on $\wt\sC_2$.
Let $X=(X_n, n \ge 0, Q^x_\om, x \in \sG(\om))$ be the
simple random walk on $\sG(\om)$. We define 
the annealed law $\bP^*$ by the semi-direct product 
$\bP^*= \bP \times Q^0_{\omega}$, and 
the rescaled height process $Z^{(n)}$ by
$$ Z^{(n)}_t = n^{-1/3} d(0,X_{\lfloor nt\rfloor }), \q t\ge 0, $$
where $d(.,.)$ is the graph distance in $\sG(\om)$.

\ms
The following summarizes the main results in of [Ke2] in the tree
case.

\proclaim Theorem 1.1. (a) ((1.19) in [Ke2].)
Let $T_N=\min\{ n: d(0,X_n)=N\}$. Then
for all $\eps>0$ there exist $\lam_1, \lam_2$ such that
$$ \bP^*( \lam_1 \le N^{-3} T_N \le \lam_2) \ge 1-\eps, \q \hbox{ for all }
 N \ge 1. $$
(b) ((1.16) in [Ke2], full proof in [Ke3].)
Under $\bP^*$ the processes $Z^{(n)}$ converges weakly in $C[0,\infty)$ to
a process $Z$ which is not the zero process.

\ms To understand why the $n^{-1/3}$ scaling arises in (b) it is helpful
to consider the behaviour of random walks on regular deterministic
graphs with a large scale fractal structure -- see for example
[Jo], [BB2], [HK], [GT1], [GT2] and [BCK]. 
Let $d_f\ge 1$ give the volume growth,
so that $|B(x,r)| \sim r^{d_f}$, and suppose that the effective electrical
resistance $R(x, B(x,r)^c)$ between $x$ and the exterior of 
$B(x,r)$ satisfies  $R(x, B(x,r)^c) \sim r^\zeta$, where $\zeta>0$. 
In this `strongly recurrent' case (see [BCK] for simple 
recent proofs using ideas
that are also used in this paper) one finds that the mean time for
$X$ to escape from $B(x,r)$ scales as $r^{d_w}$ where 
$d_w=d_f+\zeta$. 
While the IIC $\sG$ is more irregular than the sets considered in these
papers, it still has properties similar to regular graphs with 
$d_f=2$.
Further, by Proposition 2.10 below, 
only $O(1)$ points on  $\pd B(x,r/4)$ are connected to $B(x,r)^c$
by a path outside $B(x,r/4)^c$, so one has  $R(x, B(x,r)^c) \sim r$,
giving $\zeta=1$ and $d_w=3$. 

In this paper we study the simple random walk on $\sG$, and in particular
investigate both quenched and annealed properties of
its transition densities. For technical convenience we
work with the continuous time simple random walk
on $\sG$, which we denote $Y=(Y_t, t\in [0,\infty), P^x_\om, x \in \sG(\om))$.
Since we consider the law of $Y$ with general starting points $x$, we
need to consider the measures $\bP_x=\bP(\cdot| x\in \sG)$
and  $\bP_{x,y}=\bP(\cdot| x,y\in \sG)$.

Unlike [Ke2] we restrict our attention to branching processes 
with a Binomial offspring distribution. Our main reason for this
is to maintain good uniform control of the laws $\bP_x$. It
is clear by symmetry that $\bP_x( |B(x,r)| > \lam)$ is the same 
for any $x \in \bB_N$, and in fact we have uniform bounds for all
$x \in \bB$. (These probabilities are not equal for all $x$, since
a higher level $x$ is likely to be further from the backbone of
the cluster). For a general branching process, the labels
of the point $x$ may give a substantial amount of information 
about the size of the cluster near $x$.

\proclaim Theorem 1.2. (a)
There exist $c_0,c_1,c_2$, $S(x)$ such that for each $x$,
$$ \tbP_x( S(x) \ge m) \le c_0 (\log m)^{-1}, \eqno(1.2)$$
and on $\{\om: x \in \sG(\om)\}$ 
$$  c_1t^{-2/3}(\log\log t)^{-17} \le q^{\omega}_t(x,x)\le 
 c_2 t^{-2/3}(\log\log t)^{3} \hbox{ for all }
t \ge S(x).  \eqno(1.3) $$
(b) $d_s(\sG)=4/3$ $\bP$--a.s.

\smallskip The cluster $\sG$ contains large scale fluctuations, so that
$q_t(x,x)$ does have oscillations of order $(\log\log t)^c$ as
$t \to \infty$ -- see Lemma 5.1. 

\proclaim Theorem 1.3. (a) We have
$$ c_1 t^{1/3} \le \tbE_x E^x_\om d(x,Y_t) 
\le \tbE_x E^x_\om \sup_{0\le s\le t} d(x,Y_s) 
\le  c_2 t^{1/3}. \eqno(1.4) $$
(b) There exists $T(x)$ with $\tbP_x(T(x)<\infty)=1$ such that
$$ c_3  t^{1/3} (\log\log t)^{-12} \le 
   E^x_\om [d(x,Y_t)]\le c_4 t^{1/3}\log t \q \hbox{ for all }
 t\ge T(x) . \eqno(1.5) $$

\ms We also have (annealed) off-diagonal bounds for 
$q_t^\om(x,y)$. These are of the same form as the bounds
$$ ct^{-d_f/d_w} \exp( - c' (d(x,y)^{d_w}/t)^{1/(d_w-1)})$$
obtained for regular fractal graphs.

\proclaim Theorem 1.4. (a) Let $x,y \in \bB$. Then
$$ \tbE_{x,y} q^{\omega}_t(x,y) \le c_1t^{-2/3}
\exp\big( - c_2 (\frac{d(x,y)^3}t)^{1/2} \big). \eqno(1.6) $$
(b) Let $x,y \in \bB$, with $d(x,y)=R$, and $c_3 R \le t$. Then
$$ \tbE_{x,y} q^{\omega}_t(x,y) \ge c_4 t^{-2/3}
 \exp( - c_5  (R^3/t)^{1/2} ). \eqno(1.7) $$

\ms Define the continuous time rescaled height process
$$ \wt Z^{(n)}_t = n^{-1/3} d(0,Y_{nt}), \q t\ge 0. $$
By Theorem 1.3(a) the processes $(\wt Z^{(n)}, n \ge 1)$ are tight
with respect to the annealed law given by the
semi-direct product $\bP^*= \bP \times P^0_{\omega}$.
(This is much easier to prove than the full convergence given in
Theorem 1.1(b).)
However, the large scale fluctuations in $\sG$ mean
that we do not have quenched tightness.

\proclaim Theorem 1.5. $\bP$-a.s., the processes
 $(\wt Z^{(n)}, n \ge 1)$ are not tight with respect to $P^0_{\omega}$.

\ms In Section 2 we recall various properties of branching processes,
and obtain the geometrical properties of $\sG$ that we will require.
In particular we
show that, with high probability, balls $B(x,r) \subset \sG$ have roughly 
$r^2$ points, and $O(1)$ disjoint paths between $B(x,r/4)$ and
$B(x,r)^c$. Based on this, we define various types of possible `good'
behaviour of a ball $B(x,r)$, and the cluster in a neighbourhood
of the path between points $x,y \in \sG$.
In Section 3 we review some general properties of random walks on
graphs.
Our main estimates are given in Section 4, for the random walk
on a deterministic subset $\sG$ of $\bB$ for which balls and paths
are `good' in the ways given in Section 2.
 Finally, in Section 5 we tie together the results
of Sections 2 and 4, and prove 
Theorems 1.2--1.5.
%Theorems 1.2--1.4.

\ms
Throughout this article, $f_n\sim g_n$ means that 
$\lim_{n\to\infty}f_n/g_n=1$.
We use $c$, $c'$ and $c''$ to denote strictly
positive finite constants whose values are not significant and may
change from line to line. We write $c_i$ for positive constants whose values are
fixed within each theorem, lemma etc. 
When we cite a constant $c_1$ 
in Lemma $2.2$, say, we denote it as $c_{2.2.1}$.  
None of these constants depend on the degree $n_0$ of the tree.

\bigbreak \noindent
{\bf 2. The incipient infinite cluster } 

\ms 
We begin with some estimates for the critical  
Bienaym\'e-Galton-Watson branching processes 
$X_n$, $n\ge 0$, with $X_0=1$ and 
offspring distribution $\Bin(n_0,1/n_0)$ where
$n_0\ge 2$. These are quite well known, but
as we did not find them anywhere in exactly the form we needed, we
give the proofs (which are quite short) here.

Let $f$ be the generator of the
offspring distribution, so that 
$$ f(s) = E(s^{X_1}) = n_0^{-n_0}(s+n_0-1)^{n_0}. \eqno (2.1) $$
From [Har] p. 21 we have
$$ P(X_n>0) \sim \frac{2}{nf''(1)} = \frac{2n_0}{(n_0-1)n}. \eqno (2.2)$$
Let
$$ Y_n = \sum_{k=0}^n X_k, \qq g_n(s) =  E(s^{Y_n}), \qq 
 f_n(s) = Es^{X_n}. $$
Then conditioning on $X_1$ we obtain 
that $f_{n+1}(s)= f(f_n(s))$, and  
$$ g_{n+1}(s)= s f(g_n(s)) = \frac s{n_0^{n_0}}(g_n(s)+n_0-1)^{n_0}.  $$
Set 
$$ h_n(\th) = \log g_n(e^\th), \qq k_n(\th) = \log f_n(e^\th). $$

\proclaim Lemma 2.1. 
(a) Let  $1<\al \le 2$. Then
$$ h_n(\th) \le (1+\al n) \th, \q \hbox{ provided } 
 0\le \th \le  \frac{\al-1}{(1+\al n)^2}. \eqno (2.3) $$
(b) 
$$  k_n(\th) \le \th + 2n \th^2 , \q \hbox{ provided } 
 0<\th \le \frac{1}{6n}. \eqno(2.4)$$

\proof 
Note that $h_n$ and $k_n$ are continuous, strictly increasing and 
$h_n(0)=k_n(0)=0$.

For (a) we have
$$ h_{n+1}(\th)= 
\log\Big( \fract{e^\th}{n_0^{n_0}}(e^{h_n(\th)}+n_0-1)^{n_0} \Big)
 = \th + n_0 \log \frac 1{n_0}(e^{h_n(\th)}+n_0-1). $$
Let $a_n=\min\{\th: h_n(\th)=1\}$. Then since $e^x\le 1+x+x^2$
on $[0,1]$, on $[0,a_n]$,
$$   h_{n+1}(\th) \le \th + n_0 \log( 1+ \frac1{n_0} h_n(\th) 
+ \frac 1{n_0} h_n(\th)^2)
 \le \th + h_n(\th) + h_n(\th)^2.  \eqno (2.5) $$ 
We verify (2.3) by induction.
Since $h_0(\th)=\th$, (2.3) holds for $n=0$. 
Writing $b_n(\al)= (\al-1)/(1+\al n)^2$, we have
$h_n(\th)\le 1$ for  $\th \in [0,b_n(\al)]$. 
So, using (2.5) and (2.3) for $n$
$$ h_{n+1}(\th) \le  (1+\al(n+1))\th + (1+\al n)^2 \th^2 - (\al-1)\th 
\le  (1+\al(n+1))\th, $$
proving (2.3) for $n+1$.

\med (b) 
Similarly, provided $k_n(\th)\le 1$,
$$  k_{n+1}(\th) = n_0 \log\Big( 1 + \frac{e^{k_n(\th)}-1}{n_0} \Big) 
\le  k_n(\th) + k_n(\th)^2.  \eqno (2.6) $$ 
Using (2.4) for $n$ we obtain, since $\th + 2n \th^2 \le 4\th/3$,
$$ k_{n+1}(\th) \le (\th + 2n\th^2) +  (\th + 2n\th^2)^2
 \le  (\th + 2n\th^2) + 16 \th^2/9 \le  (\th + 2(n+1)\th^2), $$
proving (2.4) for $n+1$. \qed

\small {\bf Notation.} Let $\xi$ be a random variable. We write
$\lam \xi[n]$ for a r.v. with the distribution of $\lam \sum_1^n \xi_i$, where
$\xi_i$ are i.i.d. with $\xi_i \eqd \xi$. We also write $\Ber(p)$
and $\Bin(n,p)$ for the Bernoulli and Binomial distributions respectively.
Using this notation
we have for example $(\xi[n])[m]=\xi[nm]$, and $\Bin(n,p) \eqd \Ber(p)[n].$
 We write $\succcurlyeq$ for stochastic domination.

\proclaim Lemma 2.2. For any $\lambda>0$
$$  P( X_n[n] \ge \lam n ) \le c_1 e^{ -\lam/6},
\eqno(2.7)$$
$$ P( Y_n[n] \ge \lam n^2 ) \le c_2e^{ -\lam/5}.
\eqno(2.8)  $$

\proof 
Let $\th=1/6n$. Using (2.4)
$$ \eqalign{
 \log P( X_n[n] \ge \lam n) &\le -\th\lam n + n k_n(\th) \cr
 &\le -n \th(\lam-2) = -(\lam-2)/6, }$$
proving (2.7). 

Let 
If $\th\le b_n(\al)$ then 
$$ \eqalign{
  P(Y_n[n] \ge \lam n^2 ) &= P(e^{\th Y_n[n]} \ge e^{ \th \lam n^2}) 
 \le  e^{ -\th\lam n^2} Ee^{\th Y_n[n]} \cr
 &= \exp(  -\th\lam n^2 + n h_n(\th))
 \le  \exp (   -\th\lam n^2  + (1+2 n) n \th). }$$
So taking $\al=2$ and $\th=b_n(2)=(1+ 2n)^{-2}$
$$ \log  P(Y_n[n] \ge \lam n^2 ) \le 
- \frac{n^2(\lam-2)}{(1+2 n)^2} + \frac{n}{(1+2 n)^2}
 \sim -\fract15 \lam + c_3. $$
\qed

\proclaim Lemma 2.3. (a) There exist $c_0>0$, $p_0>0$ such that
$$ P(Y_n > c_0 n^2) \ge \frac{p_0}{n}. $$
(b)  If $\eta_n  \eqd \Bin(n, p_0/n)$ 
then $Y_n[n] \succcurlyeq c_0 n^2 \eta_n$.

\proof (a) This should be in literature, but is also easy to prove directly.
Let $A_n = \{ X_{n/2}>0\}$, and $a_n = P(A_n)$. Then by (2.2)
$a_n \sim (2n_0/(n_0-1))n^{-1}$. We have
$EY_n = n+1$ and $EY_n^2 \le c_1 n^3$, where $c_1$ does not depend
on $n_0$. 
On $A^c$ we have $Y_{n/2}=Y_n$, so 
$$ n+1 = EY_n = E(Y_n;A_n)+E(Y_n;A_n^c) \le E(Y_n|A_n)P(A_n) + EY_{n/2}. $$
It follows that
$$  E(Y_n|A_n) \ge \frac{n/2}{a_n} \ge c_2 n^2. $$
Also,
$$ E(Y_n^2|A_n) \le P(A_n)^{-1} E(Y_n^2;A_n) \le c_3 n^4. $$
Using the `Backwards Chebyshev' inequality 
$P(\xi \ge \half E\xi) \ge (E \xi)^2/(4E \xi^2)$ 
with respect to $P(\cdot|A_n)$ then gives 
$$   P(Y_n> \half c_2 n^2 |A_n)\ge  P(Y_n> \half E(Y_n|A_n)|A_n) \ge 
\frac{c_2^2 n^4}{4 c_3 n^4} = c_4. $$
So $P(Y_n> \half c_2 n^2) \ge  P(Y_n> c_2n^2 |A_n) P(A_n) 
\ge c_4 a_n \ge c_5 n^{-1}$, and taking $c_0=\half c_2$,
$p_0=c_5$, this proves (a).
 
\noi (b) Let now $Y^{(j)}_n$ be i.i.d. copies of $Y_n$, and
$F_j=\{  Y^{(j)}_{n}> c_0 n^2 \}$. Then if 
$\xi_j = 1_{F_j}$,
by (a) we have $P(\xi_j=1) \ge p_0/n$. So, 
$$ Y_n[n] = \sum_{j=1}^n Y^{(j)}_n 
\succcurlyeq \sum_{j=1}^n c_0 n^2 \xi_j  \succcurlyeq
 c_0 n^2 \eta_n, $$
proving (b). \qed

\proclaim Lemma 2.4. For $0<\lam<1$,   
$$  \exp( - c_1/\lam) \le 
P( Y_n[n] \le \lam n^2) \le \exp( - c_2/\lam^{1/2}). \eqno(2.9)$$

\proof To prove the upper bound let $c_0=c_{2.3.0}$, and 
$m=(\lam/c_0)^{1/2}n$. Using Lemma 2.3 we have 
$$ Y_n[n] = \sum_{i=1}^n Y^{(i)}_m 
\succcurlyeq \sum_{i=1}^n c_0 m^2 \xi_i = \lam n^2  \sum_{i=1}^n \xi_i;$$
here $\xi_i$ are i.i.d. ${\tt Ber}(p_0/m)$ r.v.
So  
$$ P(Y_n[n] < \lam n^2) \le P( \sum_{i=1}^n \xi_i  <1)
 = (1- p_0/m)^n \le \exp( -p_0 n/m ) =  \exp( - c_0^{1/2} p_0/\lam^{1/2}). $$

For the lower bound let $k\ge 1$ and $m=n/k$. 
Let  $G_j=\{  X^{(j)}_{m}=0\}$,
and $G=\cap_{1\le j\le n} G_j$. Then $P(G_j)\ge (1-c/m)^n$ so
$$ \eqalign{
 P(Y_n[n] < \lam n^2) &\ge   P(Y_n[n] < \lam n^2|G)P(G)  \cr
 &\ge (1-c/m)^n \Big( 1 -  P(Y_n[n] > \lam n^2|G) \Big) \cr
 &\ge c'e^{-c''k} \Big( 1 -  P(Y_n[n] > \lam n^2|G) \Big) . }$$
On $G$ we have $Y_n[n] =\sum_{j=1}^n Y^{(j)}_m$, so
$$  \eqalign{
 P(Y_n[n] > \lam n^2|G) \le 
   \frac{ E(\sum_{j=1}^n Y^{(j)}_m |G) }{ \lam n^2} 
=  \frac{ n E (Y^{(1)}_m |G_1) }{ \lam n^2} 
 \le  \frac{ E Y^{(1)}_m }{\lam n P(G_1)}\le  \frac{c}{k\lam}. }$$
Taking $k$ such that $c/(k\lam)=\half$ completes the proof.  \qed

\ms We will need to consider the following modified branching process.
Let $\wX=(\wX_n, n\ge 0)$ be a branching process with $\wX_0=1$ and the same
${\Bin}(n_0, 1/n_0)$ offspring distribution as $X$, except that
at the first generation we have $\wX_1 \eqd {\Bin}(n_0-1, 1/n_0)$.

\proclaim Lemma 2.5. (a)  For any $\lambda>0$
$$  P( \wX_n[n] \ge \lam n ) \le c_1 e^{ -c_2 \lam},
\eqno(2.10)$$
$$ P( \wY_n[n] \ge \lam n^2 ) \le c_3 e^{ -c_4\lam}.
\eqno(2.11)  $$
(b)  For $0<\lam<1$,   
$$  \exp( - c_5/\lam) \le 
P( \wY_n[n] \le \lam n^2) \le \exp( - c_6/\lam^{1/2}). \eqno(2.12)$$
(c)  There exists $p_1>0$ such that 
$\wY_n[n] \succcurlyeq c_7 n^2 \Bin(n,p_1/n)$.

\proof (a) and the 
lower bound in (b) are
immediate from Lemmas 2.2 and 2.4, since $\wX_n \preccurlyeq X_n$ and 
$\wY_n \preccurlyeq Y_n$.

 For the upper bound in (b), we can write
$$ \wY_n[n]= n + \sum_{i=1}^M Y^{(i)}_{n-1}, $$
where $M \eqd \Bin(n(n_0-1), 1/n_0)$, and 
$Y^{(i)}$ are independent copies of $Y$. Similarly,
$$ Y_m[m]= m + \sum_{i=1}^{M'} Y^{(i)}_{m-1}, $$
where $M' \eqd \Bin(nn_0, 1/n_0)$.
So if $m=n (n_0-1)/n_0$ then 
$$  \wY_n[n]= n + \sum_{i=1}^M Y^{(i)}_{n-1} 
 \ge m +  \sum_{i=1}^M Y^{(i)}_{m-1} = Y_m[m]. \eqno(2.13)$$
(2.12) now follows from Lemma 2.4, since $\half n \le m \le n$.
\nl (c) We have $ \Ber(p)  \succcurlyeq \half \Ber(p/2)[2]$. 
So, using (2.13), with $m$ as in (b),
$$  \eqalign{
\wY_n[n] \succcurlyeq  Y_m[m] 
&\succcurlyeq  c_0 m^2 \Bin(m, p_0/m)  \cr
&\succcurlyeq  \half c_0 m^2 \Bin(2m, p_0/2m) \cr
&\succcurlyeq  \half c_0 m^2 \Bin(n, p_0/2m)
 \succcurlyeq c_1 n^2 \Bin(n, p_1/n) .} $$
\qed

\bs We now define the random graph $\sG$ we will be working with.  We
could regard this either as critical percolation on the $n_0$-ary tree
$\bB$, conditioned on the cluster containing the root 0 being
infinite, or as the (critical) Bienaym\'e-Galton-Watson process with
${\Bin}(n_0,1/n_0)$ offspring distribution, conditioned on
non-extinction.

Let $\bB$ be the $n_0$-ary tree, and let $0$ be the root. 
A point $x$ in the $n$th generation (or level) is written
$x=(0,l_1, \cdots, l_n)$, where $l_i\in \{1,2,\cdots, n_0\}$. Let $\bB_n$
be the set of $n_0^n$ points in the  $n$th generation, and let
$\bB_{\le n}=\cup_{i=0}^n \bB_i$. If $x\in \bB_k$ we write $|x|=k$.
If $x =(0,l_1, \cdots, l_n) \in \bB_n$, let 
$a(x,r)=(0,l_1, \cdots, l_{n-r})$ be the ancestor of $x$ at level
$|x|-r$. 

We regard $\bB$ as a graph (in fact a tree) with edge set
$E(\bB) =\big\{ \{x,a(x,1)\}, x\in \bB-\{0\} \big\}$. 
Let $\eta_e$, $e\in E(\bB)$, be i.i.d. Bernoulli $1/n_0$ r.v. defined on a
probability space 
$(\Omega, \sF, P)$. 
If $\eta_e=1$ we say the edge $e$ is {\sl open}. Let
$$ \sC(0) = \{x \in \bB: \hbox{ there exists an $\eta$--open path from 
 $0$ to $x$} \} $$
be the open cluster containing $0$. 
It is clear that $Z_n = | \sC(0) \cap \bB_n|$ is a critical 
GW process with $\Bin(n_0,1/n_0)$ offspring distribution. 
Here and in the following, $|A|$ is a cardinality of the set $A$. 
As $Z$ has extinction probability 1, the cluster $\sC(0)$ is 
$P$--a.s. finite. 

We have

\proclaim Lemma 2.6. ([Ke2,  Lemma 1.14]). 
Let $A\subset \bB_{\le k}$. Then 
$$ 
\lim_{n\to\infty} P(\sC(0)\cap \bB_{\le k} =A|Z_n \ne 0)
 = |A \cap \bB_k| P(\sC(0)\cap \bB_{\le k} =A), \eqno(2.14)$$
and writing  
$\bP_0(A)= |A \cap \bB_k| P(\sC_{\le k} =A)$,
$\bP_0$ has a unique extension to a probability measure $\bP$ 
on the set of infinite connected subsets of $\bB$ containing $0$.

Let $\sG'$ be a rooted labeled tree chosen with the distribution
$\bP$: we call this the {\sl incipient infinite cluster} (IIC) on $\bB$.
For more information on $\sG'$ see [Ke2] and [vH] but we remark that
$\bP$--a.s. $\sG'$ has exactly one infinite descending path from
$0$, which we call the {\sl backbone}, and denote $H$. 

\smallskip 
It will be useful to give another construction of the IIC, 
obtained by modifying the cluster $\sC(0)$ rather than its law.
We can suppose the probability space 
$(\Omega, \sF, P)$
carries i.i.d.r.v.  $\xi_i$, $i\ge 1$ uniformly distributed on
$\{1,2,\cdots, n_0\}$, and independent of $(\eta_e)$. For $n\ge 0$ let
$\Xi_n = (0, \xi_1, \dots , \xi_n)$, and let
$$ \wt \eta_e = \cases{ 1 & if $e=\{\Xi_n, \Xi_{n+1} \}$ for
     some $n\ge 0$, \cr
     \eta_e & otherwise. \cr}$$
Then (see [vH]) if 
$$ \sG = \{x \in \bB: \hbox{ there exists a $\wt\eta$--open path from 
 $0$ to $x$} \}, $$
$\sG$ has law $\tbP$. It is clear that the backbone
of  $\sG$ is the set $H=\{ \Xi_n, n\ge 0\}$.  

\ms For $x, y \in \bB$ let  
$$  \tbP_x(\cdot)=\tbP (\cdot| x\in \sG), \qq
 \tbP_{xy}(\cdot)=\tbP (\cdot| x,y \in \sG), $$
and let $\tbE_x$ and  $\tbE_{xy}$ denote expectation
with respect to  $\tbP_x$ and  $\tbP_{xy}$ respectively. 
Given a descending path  $b=\{0,b_1, b_2,\dots\}$, 
(which we call a {\sl possible backbone}) let
$$ \tbP_{x,b}(\cdot)  = \tbP(\cdot | x \in \sG, H=b), $$
and define $\tbP_{x,y,b}$ analogously. 

For each $x$, $y\in \bB$, let $\gam(x,y)$ be the unique geodesic path
connecting $x$ and $y$. 
We say that $z$ is a {\sl middle point} of $\gam(x,y)$ if
$z\in \gam(x,y)$ and $|d(x,z)-\half d(x,y)|\le \half$.
We remark that the construction of $\sG$ makes
it clear that $\tbP_{x,y,b}(\eta_e=1)=1$ if the edge $e$ lies in any
of the paths $b$, $\gamma(0,x)$ and $\gamma(0,y)$, and that under
$\bP_{x,y,b}$ the r.v. $\eta_e$, $e \not\in b \cup \gamma(0,x) \cup
\gamma(0,y)$ are i.i.d. with $\bP_{x,y,b}(\eta_e=1)=1/n_0$.

\med {\bf Notation.} 
We consider the tree $\sG=\sG(\omega)$.  
Let $d(x,y)$ be the graph distance between $x$ and $y$, and
$$ B(x,r) =\{y\in \sG : d(x,y) \le r\}. $$
We write $D(x)$ for the set of descendants of $x$. More precisely,
$y\in D(x)$ if and only if  $x\in \gam(0,y)$. Note that $x\in D(x)$.
If $y\in D(x)$ we call $x$ an {\sl ancestor} of $y$ and $y$
a {\sl decedent} of $x$.
We set
$$ D_r(x) = \{ y \in D(x): d(x,y) = r \}, \qq
 D_{\le r}(x) = \cup_{i=0}^r D_i(x).  $$
We also set
$$ D(x;z) = \{ y \in D(x): \gamma(x,y)\cap \gamma(x,z)=\{x\} \},$$
and write $D_r(x;z)=D_r(x)\cap D(x;z)$, 
 $D_{\le r}(x;z) = D_{\le r}(x)\cap D(x;z)$.
Thus if $z\in D(x)$ then $y\in D(x;z)$ if and only if the
lines of descent from $x$ to $y$ and $z$ are disjoint, except
for $x$. (Note that $D(x;x)=D(x)$.)
For any $A\subset \sG$ we write 
$$ \pd A= \{ y\in \sG-A: y\sim x \hbox{ for some } x\in A\}. $$

\ms The estimates at the beginning of this Section lead to volume 
growth estimates for $\sG$. For $x\in \sG$ let $\mu_x$ be the degree
of $x$, and for $A\subset \sG$ set $\mu(A)=\sum_{x \in A} \mu_x$.
We write
$$ V(x,r) =\mu(B(x,r)). $$
Note that as $\sG$ is a tree, we have
$$ |B(x,r)| \le V(x,r) \le 2|B(x,r+1)|. \eqno (2.15) $$

\proclaim Proposition 2.7. 
(a)  Let $\lambda>0$, $r\ge 1$ and $x,y \in \bB$, and $b$ be a
possible backbone.  Then
$$\tbP_{x,y,b} ( V(x,r) > \lambda r^2)\le c_0 \exp (-c_1 \lambda), 
\eqno(2.16)$$
and 
$$\tbP_{x,y,b} (V(x,r)<\lambda r^2)\le c_2 \exp (-c_3/\sqrt{\lambda}).
\eqno(2.17)$$
\nl (b) The bounds (2.16) and (2.17) also hold for the laws
$\tbP_{x,b}$, $\tbP_{x,y}$, and $\tbP_{x}$. 

\proof It is enough to prove (a), since the bounds for
$\tbP_{x,b}$ follow by taking $y=0$, and those for 
$\tbP_{x,y}$ and $\tbP_{x}$ then follow on integrating over $b$. 
Also, using (2.15), it is enough to bound $|B(x,r)|$.

We will assume that $|x|>r$; if not we can use the same arguments
with minor modifications. 
Let $x_i= a(x,i)$ for $0\le i \le r$.
If the backbone intersects $B(x,r)$ then let $s$ be the smallest
$i$ such that $x_i\in H$, and let $v_0=x_s$ and $v_i$, $i\ge 1$
be the backbone descending from the point $v_0$. Similarly
if $\gamma(0,y)$  intersects $B(x,r)$ then let $t$ be the smallest
$j$ such that $y_j\in \gamma(0,y)$, and let  
$w_0=y_t$ and $w_i$, $1\le i \le t$ be the path $\gamma(w_0,y)$.

Then we have
$$ B(x,r) \subset \big(\cup_{i=0}^{r} D_{\le r}(x_i;x) \big)\cup
 \big(\cup_{i=1}^r D_{\le r}(v_i;v_{3r}) \big)
  \cup \big(\cup_{i=1}^{r\wedge t} D_{\le r}(w_i;y) \big).$$
Under $\tbP_{x,y,b}$ the r.v.  $|D_{\le r}(\cdot;\cdot)|$ above
are i.i.d., with the same law as $\wY_r$. 
Thus $|B(x,r)| \preccurlyeq \wY_r[r][3]$, and 
by Lemma 2.5(a),
$$  \tbP_{x,y,b}( |B(x,r)| > \lambda r^2)\le c \exp (-c' \lambda). $$

The proof of (2.17) is very similar. We have
$ \cup_{i=0}^{r/2} D_{\le r/2}(x_i;x) \subset B(x,r), $
so that $|B(x,r)| \succcurlyeq \wY_{r/2}[r/2]$, and using
Lemma 2.5(b) leads to (2.17). \qed

\ms We also wish to show that oscillations in 
$n^{-2}V(0,n)$ exist. If $W \eqd {\Bin}(n,p/n)$ then
straightforward calculations give that
$$ P(W=k) \ge c_0 e^{-k \log(k/p)}, \q 0\le k \le n^{1/2}. \eqno(2.18)$$

\proclaim Proposition 2.8.
(a) For any $\eps >0$  
$$  \limsup_{n \to \infty} \frac{ V(0,n)}{ n^2(\log\log n)^{1-\eps}}=\infty, 
 \q \bP-a.s. $$
(b) There exists $c_0<\infty$ such that
$$ \liminf_{n \to \infty} \frac{(\log\log n) V(0,n)}{n^2} 
\le c_0.  \q \bP-a.s. $$

\proof It is enough to prove these for the law $\bP_b$, for
any fixed possible backbone $b=\{0, y_1, y_2, \dots \}$. 
\nl (a)
Let
$$ Z_n = |\{x: x\in D(y_i;y_{i+1}),\, d(x,y_i)\le 2^{n-2}, \,
2^{n-1}\le i \le 2^{n-1} + 2^{n-2} \}|. $$
Thus $Z_n$ is the number of descendants off the backbone, 
to level $2^{n-2}$, of points $y$ on the backbone between levels $2^{n-1}$ and 
$2^{n-1} + 2^{n-2}$. So $|B(0,2^n)|\ge Z_n$, the r.v.
$Z_n$ are independent, and $Z_n \eqd \wY_{2^{n-2}}[2^{n-2}]$.
Using Lemma 2.5(c) we have, if $a_n= (\log n)^{1-\eps}$,
and $\eta_n \eqd \Bin(n,p_1/n)$,
$$ \eqalign{
\bP_b(|B(0,2^n)|\ge a_n4^n) &\ge  \bP_b( Z_n\ge a_n4^n) \cr
 &\ge P( \wY_{2^{n-2}}[2^{n-2}] \ge  a_n4^n) \cr
 &\ge P( \eta_{2^{n-2}} \ge a_n) \ge c e^{- a_n \log a_n}. }$$
As $Z_n$ are independent, (a) follows by the second Borel-Cantelli
Lemma. 

\nl (b) Let $n_k = \exp( 2 k \log k)$, so that $k^2 n_{k-1}\le n_k$,
and let
$$ W_k = \bigcup_{i=0}^{n_k-1} D(y_i;y_{n_k}), \q
 V_k = D_{\le n_k-n_{k-1}}(y_{n_{k-1}}). $$
Then the r.v. $|V_k|$ are independent and 
$ B(0,n_k) \subset W_{k-1} \cup V_k$.

Fix $0<\ep<1/3$ and let 
$$F(i,k)= \{ D_{k^{1+\ep} n_k}(y_i;y_{i+1})=\emptyset\}. $$
Then since $X_n \succcurlyeq  \wX_n$
$$ \bP(F(i,k))= P( \wX_{k^{1+\eps} n_k }=0) \ge  P( X_{k^{1+\eps} n_k }=0)
\ge 1 - \frac{c}{k^{1+\eps} n_k}. $$
Let $G_k= \cap_{i=0}^{n_k-1} F(i,k)$; we have
$$ \bP(G_k^c) \le c/k^{1+\ep}. $$
On the event $G_k$ we have that $|W_k|$ is stochastically
dominated by $\sum_{i=1}^{n_k} Y^{(i)}_{k^{1+\ep} n_k}$, so
$$ \eqalign{
 \bP(|W_k| \ge  k^3 n_k^2) &\le  \bP(G_k^c) 
  + P( Y_{k^{1+\ep} n_k}[ k^{1+\ep} n_k] \ge k^{1-2\ep} (k^{1+\ep} n_k)^2)\cr
 &\le c k^{-2} + e^{-c'k^{1-2\ep}} \le c''k^{-2}. }$$
Thus $|W_k|\le k^3 n_k^2$ for all large $k$.
Now $|V_k| \preccurlyeq Y_{n_k}[n_k]$, so
$$ \bP( |V_k| < c_1 (\log k)^{-1} n_k^2)
 \ge P( Y_{n_k}[n_k] <  c_1 (\log k)^{-1} n_k^2) \ge e^{-c \log k}
\ge k^{-1} $$
if $c_1$ is chosen large enough. As the r.v. $|V_k|$ are independent,
we deduce that $|V_k|<  c_1 (\log k)^{-1} n_k^2$ for all $k$ in an infinite
set $J$. For all large $k\in J$,
$$ |B(0,n_k)| \le |V_k|+ (k-1)^3 n_{k-1}^2
  \le  (c_1 (\log k)^{-1} + k^{-1}) n_k^2
\le \frac{2 c_1 n_k^2}{\log\log n_k} . $$
\qed

\med {\bf Remark.} Let $\sC_\infty$ denote the unique infinite cluster
for supercritical bond percolation (i.e. $p>p_c$) in $\bZ^d$. 
Then writing $Q(x,N)$ for the box side $N$ and center $x$
$$ \frac{ | \sC_\infty \cap Q(x,N)|}{|Q(x,N)|} \to \th(p). $$
Propositions 2.7 and 2.8 show that one does not get this kind
of convergence for $\sG$, which is a much more irregular set than
the clusters considered in [B2].

\med {\bf Definition 2.9.} Let $x\in \sG$, $r\ge 1$. 
Let $M(x,r)$ be the smallest number $m$ such that there exists
a set $A=\{z_1, \dots, z_m\}$ with $d(x,z_i)\in [r/4,3r/4]$ for
each $i$, such that any path $\gam$ from $x$ to $B(x,r)^c$
must pass through the set $A$. (Since $\sG$ is a tree, the best
choice of such a set $A$ will in fact have the points at a distance
$r/4$ from $x$, but we will not need this.)

\proclaim Proposition 2.10.
There exist $c_1,c_2>0$ such that for each $r\ge 1$ and each $x, y\in \bB$,
and possible backbone $b$
$$ \tbP_{x,y,b} (M(x,r)\ge m)\le c_1e^{-c_2 m} .$$
Similar bounds hold for $\tbP_{x,y}$, $\tbP_{x,b}$ and $\tbP_x$.

\proof We just consider the case $y=0$; the general case is similar
but a little more complicated since we would also need to consider offspring
on the branch $\gamma(0,y)$.
Let $w_0=a(x,r/3)$. If $w_0\in b$ then let $w_1$ be the point in the backbone 
at level $|x|+r/3$, otherwise let $w_1=w_0$. 
Let 
$$ A_1 = \cup_{z \in \gam(w_0,x), z\not\in b} D_{r/4}(z;x),
\qq A_2 = \cup_{z \in \gam(w_0,w_1), z\neq w_1} D_{r/4}(z;w_1). $$
Let $N_i=|A_i|$; we have $N_1  \preccurlyeq X_{r/4}[1+r/4]$ and
$N_2  \preccurlyeq X_{r/4}[r/2]$.
Now let
$$ A_i^* = \{ z \in A_i: D_{r/4}(z) \neq \emptyset\}. $$
Then any path from $x$ to $B(x,r)^c$ must pass through
$A_1^* \cup A_2^* \cup\{ w_0, w_1\}$, so 
$M=M(x,r) \le 2 + |A_1^*|+|A_2^*|$. 

Let $p_r= P(z \in A_i^*|z \in A_i)=P( X_{r/4}>0)$, so that
$p_r \le c/r$. So, if
$\kappa_i$ are i.i.d. $\Ber(p_r)$ r.v. independent of $N_i$, we have
$$ |A_i^*| \eqd \sum_{j=1}^{N_i} \kappa_j. $$
Let 
$$ W_n = \sum_{i=1}^n (\kappa_i-p_r); $$
then $W=\{W_n\}$ is a martingale, $W_n-W_{n-1}\le 1$,
$\langle W \rangle_n = np_r(1-p_r)$, and
$|A_i^*| \eqd W_{N_i} + N_i p_r$.
Choose $r$ large enough so that $p_r<\half$. Then 
$$  \tbP_{x,b}(  |A_i^*| \ge m) \le  \tbP_{x,b}(W_{N_i}+N_ip 
\ge m, N_i p \le m/2) 
  +  \tbP_{x,b}( N_i p> m/2). \eqno(2.19)$$
For the first term in (2.19) we have
$$ \eqalign{
 \tbP_{x,b}(W_{N_i}+N_i p \ge m, N_ip \le m/2) &\le
   \tbP_{x,b}(W_{N_i} \ge m/2,\langle W \rangle_{N_i} \le m(1-p)/2) \cr
&\le \exp(- \frac{(m/2)^2}{ 2( (m/2)+ m(1-p)/2)}) \le e^{-c m}, }$$
where we used an exponential martingale inequality -- see (1.6) in [F].
For the second term, note that $N_i \preccurlyeq (X_{r/4}[r/4])[2]$
and so using Lemma 2.2 we deduce that
$$  \tbP_{x,b}( N_i p > m/2)  \le  ce^{-c_3 m}. $$
Combining these bounds completes the proof. \qed

\med {\bf Definition 2.11.} Let $x\in \bB$, $r\ge 1$,
$\lam\ge 64$. We say that $B(x,r)$ is {\sl $\lam$--good}
if:
\nl (a) $x \in \sG$
\nl (b) $r^2\lam^{-2} \le V(x,r) \le r^2\lam$.
\nl (c) $M(x,r) \le \frac{1}{64} \lam$.
\nl (d) $V(x, r/\lam) \ge r^2 \lam^{-4}$.
\nl (e) $V(x, r/\lam^2) \ge r^2 \lam^{-6}$.

\proclaim Corollary 2.12.
For $x\in \bB$ and any possible backbone $b$
$$  \tbP_{x,b}( B(x,r) \hbox{ is not $\lam$--good}) \le c_1 e^{-c_2 \lam}. 
\eqno(2.20) $$

\proof By Propositions 2.7 and 2.10 the probability of 
each of conditions (a)--(d) above failing is bounded by
$\exp(-c \lam)$. \qed

\ms We now need to introduce some more complicated conditions on the 
tree $\sG$, and will prove that these hold with high probability.
These conditions describe various kinds of `good' behaviour of balls
with centers on a path $\gamma(x,y)$, and will be used when we
consider off-diagonal bounds on the transition probabilities of the
random walk in Sections 4 and 5.

Fix  $\lambda_1 \ge 64$ large enough so that the 
right hand side of (2.20) is less than $\frac14$.
For $x,y \in \bB$ and $k\in \bN$, define the event
$$ \eqalign{
 F_1(x,y,r,k) =\{ &x, y \in \sG \hbox{ and there exist at least $k$ 
 disjoint balls} \cr
 & B(z,r/2) \hbox{ with $z\in \gamma(x,y)$ and which
are $\lam_1$--good}.\} } $$
For $x,y\in  \bB$, let $z_0$ be a middle point of $\gamma (x,y)$.
Define the events 
$$ \eqalign{
A_*(z,r,N) &=\{z \in \sG \hbox{ and } B(z,r) \hbox{ is $N$--good}.\},\cr
F_*(x,y,R,k;r,N)&=  F_1(x,z_0,R,k/2) \cap F_1(z_0,y,R,k/2) \cr
  & \qq \cap A_*(x,r,N)\cap A_*(z_0,r,N)\cap A_*(y,r,N). }$$

\med {\bf Definition 2.13.} 
The vertex $x\in \bB$ satisfies the condition $G_2(N,R)$ if:
\nl (a) $x \in \sG$,
\nl (b) For every $z\in \pd B(x,NR)$ the event $F_1(x,z,R,\frac18 N)$
holds.

\proclaim Proposition 2.14. Let $x_0, y_0 \in \bB$, and $b$ be a possible backbone.

\noi (a) For $R\ge 1$, $N\ge 8$,
$$ \tbP_{x_0,y_0,b}\big(x_0 \hbox{ satisfies the condition $G_2(N,R)$}\big) 
\ge 1 - c_1 \exp(-c_2 N). $$
(b) The same bounds as in (a) hold for the laws
$\tbP_{x_0,b}$, $\tbP_{x_0,y_0}$, and $\tbP_{x_0}$. 
\nl (c) For $x_0,y_0 \in \bB$, $8\le N< d(x_0,y_0)/8$, $r\ge 1$,
$$ \tbP_{x_0,y_0,b}\big(F_*(x_0,y_0,\frac {d(x_0,y_0)}N, \fract18 N;r,N)\big) 
\ge 1 - c_3 \exp(-c_4 N). $$

\proof (a) We prove this for $y_0=0$; as in Proposition 2.10 the general
case is handled by a similar argument.

Let 
$$ F_0(y,s)=\{ y \in \sG \hbox{ and $B(y,s)$ is $\lam_1$--good.}\}, $$
and write $v_i= a(x,i)$, $R'=RN/4$.
We assume that $|x| \ge NR$ and $v_{R'}$ is on the backbone $b$:
the other cases can be handled by minor modifications to the arguments
below. 
Let $w_0$ be the highest level point in both $b$ and $\gamma(0,x)$,
and $w_i$, $i\ge 1$ be the backbone $b$ from $w_0$ on.

Under $\tbP_{x,b}$ the events $F_0(v_{Rj},\frac R2)$,
$1\le j \le N$ are independent, and  $\bP_{x,b}(F_0(v_{Rj},\frac R2)^c )$
$\le \frac14$.
So standard exponential bounds give 
$$ \bP_{x,b}( F_1(x, v_{R'},R,N/8)^c ) \le  c\exp(-c' N). \eqno(2.21) $$
Similarly
$$ \bP_{x,b}( F_1(w_0, w_{R'},R,N/8)^c ) \le  c\exp(-c' N). $$

Now let $A_1=\{ v_i, 0\le i \le R' \}\cup \{ w_i, 0\le i \le R' \}$;
note that under $\bP_{x,b}$ this set is non-random. 
Let  
$$ A_2= \big\{ y \in \bB: a(y,R') \in A_1, 
\gamma(y,a(y,R')) \cap A_1= \{a(y,R') \} \big\}. $$
For $y\in A_2$ let
$$  \eqalign{
H_1(y) &= F_1( a(y,R),a(y,R'),R, N/8)^c, \cr
H_2(y) &= \{ y \in \sG, \, D_{R'}(y) \ne \emptyset \}. }$$  
Then
$$ \eqalign{
  \bP_{x,b}( \bigcup_{y\in A_2} H_1(y)\cap H_2(y) ) 
 &\le \sum_{y\in A_2}  \bP_{x,y,b}( H_1(y)\cap H_2(y))
 \bP_{x,b}(y \in \sG). \cr }$$
Under $\bP_{x,y,b}$ the events $H_1(y)$ and $H_2(y)$ are independent, and
as in (2.21) we obtain  $\bP_{x,y,b}(H_1(y)) \le c\exp(-c' N)$. So,
$$  \eqalign{
\bP_{x,b}( \bigcup_{y\in A_2} H_1(y)\cap H_2(y) )
 &\le  ce^{-c' N} \sum_{y\in A_2}  
             \bP_{x,y,b} (H_2(y)) \bP_{x,b}(y \in \sG) \cr
 &=  ce^{-c' N}  \sum_{y\in A_2}  \bP_{x,b} (H_2(y)) \cr
 & = ce^{-c' N} \bE_{x,b}  \sum_{y\in A_2} 1_{H_2(y)}. }$$
The final sum above is bounded by a constant $c'$ by the same
argument as in Proposition 2.10.  

Finally, we have
$$ \eqalign{ 
  \{ \hbox{ $G_2(N,R)$ fails for $x$} \} &\subset \cr
   F_1(x, v_{R'},R,N/8)^c \; \cup  &F_1(w_0, w_{R'},R,N/8)^c \;
 \cup  \bigcup_{y\in A_2} (H_1(y)\cap H_2(y)),} $$
so combining the bounds above completes the proof. 
\noi (b) follows on integrating the bounds in (a).

For (c), we first note that, by the argument for (2.21),
$$\tbP_{x,y,b}\big(F_1(x,y,\frac {d(x,y)}N,\frac 1{16} N)^c\big) \le c' \exp(-c N).$$ 
So, using Corollary 2.12, we have 
$$ \eqalign{ 
\tbP_{x,y,b}\big(F_*^c\big) \le &\tbP_{x,y,b}\big(F_1(x,z_0,\frac {d(x,y)}N,
\frac 1{16} N)^c\big) +
\tbP_{x,y,b}\big(F_1(z_0,y,\frac {d(x,y)}N,\frac 1{16} N)^c\big)\cr
&\,+\sum_{w=x,z_0,y}\tbP_{x,y,b}\big(A_*(w,r,N)^c\big)\cr 
\le &2c' \exp(-c N)+3c' \exp(-c N)= 5c' \exp(-c N).}$$ 
\qed

\med {\bf Definition 2.15.} Let $x, y \in \bB$, $m, \th \in \bN$. 
Define the condition $G_3(x,y,m,\kappa)$ as follows. 
Let $r=d(x,y)/m$, and let $z_0=x, z_1, \dots, z_m=y$ be points on the path
$\gam(x,y)$ with $|d(z_{i-1},z_i)-r|\le 1$. 
(We choose these points in some fixed way -- for example so that 
$d(z_{i-1},z_i)$ are non-decreasing.) For each $i=1, \dots m$ let
$\Th_i$ be the smallest integer $\lam\ge \max(64, 3 c_{4.7.2}^{-1})$ 
such that $B(z_i, \lam^{20} r)$
is $\lam$--good, and $|B(z_i, r)|\ge r^2/\lam^2$.
Then  $G_3(x,y,m,\kappa)$ holds if:
\nl (a) $x,y  \in \sG$,
\nl (b) $\sum_{i=1}^m \Th_i^{54} \le \kappa m$.

\proclaim Proposition 2.16. For each backbone $b$ and $x, y \in \bB$
$$ \tbP_{x,y,b}\big( \hbox{ $G_3(x,y,m,\kappa)$ holds }\big)
\ge 1 - c_1 \kappa^{-1}. $$

\proof By Proposition 2.7 and  Corollary 2.12, 
$\tbP_{x,y,b}(\Th_i=k) \le e^{-c k}$. Thus
$\tbE_{x,y,b} \Th_i^{54} \le c'$, and so 
$$  \tbP_{x,y,b}\big( \hbox{ $G_3(x,y,m,\kappa)$ fails }\big) 
 =  \tbP_{x,y,b}\big(  \sum_{i=1}^m \Th_i^{54} > \kappa m) \le c'/\kappa. $$
\qed

\bigskip\noindent
{\bf 3. Markov chains on weighted graphs and trees}

\ms Let $\Gamma$ be a infinite connected locally finite graph.
Assume that the graph $\Gamma$ is endowed by a weight (conductance) 
$\mu_{xy}$, which is a symmetric nonnegative function on $\Gamma\times
\Gamma$ such that $\mu_{xy}>0$ if and only if $x$ and $y$ are connected 
by a bond (in which case we write $x\sim y$). We call the pair
$(\Gamma, \mu)$ a {\sl weighted graph}. We can also regard it as an 
electrical network, in which the bond $\{x,y\}$ has conductance $\mu_{xy}$.
We will be mainly concerned with the case when $\mu_{xy}=1$ if and only
if $\{x,y\}$ is an edge: we call these the {\sl natural weights}
on $\Gamma$.
Let $\mu_x=\sum_{y\in \Gamma}\mu_{xy}$ for each $x\in \Gamma$,
and set $\mu (A)=\sum_{x\in A}\mu_x$ for each $A\subset \Gamma$,
so that $\mu$ is then a measure on $\Gamma$. 

We next define a quadratic form $\sE$ on $\Gamma$ by
$$  \sE(f,g)=\frac 12\sum_{{x,y\in \Gamma}\atop
{x\sim y}}(f(x)-f(y))(g(x)-g(y))\mu_{xy},$$
and set
$$ H^2=H^2(\Gamma,\mu) =\{ f\in  \bR^{\Gamma}: \sE(f,f)<\infty\}. $$
For $f,g \in H^2$ we define $\sE(f,g)$ by polarization.
We sometimes abbreviate $\sE(f,f)$ as $\sE(f)$.
Note that if $f=\min_{1\le i\le n} g_i$ then
since 
$$ |f(x)-f(y)|^2 \le \max_i |g_i(x)-g_i(y)|^2 \le \sum_i  |g_i(x)-g_i(y)|^2, $$
it follows that
$$ \sE(f,f) \le  \sum_{i=1}^n \sE(g_i,g_i). \eqno(3.1)  $$

Let $Y=\{Y_t\}_{t\ge 0}$ be the continuous time random walk on $\Gamma$
associated with $\sE$ and the measure $\mu$. When the natural weights are 
given on $\Gamma$, $Y$ is called the simple random walk 
on $\Gam$. $Y$ is the Markov process with generator 
$$\sL f(x)=\frac 1{\mu_x}\sum_y\mu_{xy}(f(y)-f(x));$$
$Y$ waits at $x$ for an exponential mean $1$ random 
time and then moves to a neighbour $y$ of $x$ with probability
proportional to $\mu_{xy}$.
We define the transition 
density (heat kernel density) of $Y$ with respect to $\mu$ by
$$q_t(x,y)=\bP^x(Y_t=y)/\mu_y.\eqno (3.2)$$
If $A \subset \Gam$ we write
$$ T_A =\inf\{t\ge 0: Y_t \in A\}, \qq \tau_A = T_{A^c}. $$

The natural metric on the graph, obtained by counting the number of
steps in the shortest path between points, is written $d(x,y)$ for
$x,y\in \Gamma$. 
As before, we write 
$$ B(x,r)=\{y: d(x,y)\le r\}, \q V(x,r)=\mu(B(x,r)). $$

Let $A,B$ be disjoint subsets of $\Gamma$. 
The effective resistance between $A$ and $B$ is defined by:
$$R(A,B)^{-1}=\inf\{\sE(f,f): f\in H^2, f|_A=1, f|_B=0\}. \eqno(3.3) $$
Let $R(x,y)=R(\{x\},\{y\})$, and $R(x,x)=0$. In general
$R$ is a metric on $\Gam$ -- see [Kig] Section 2.3. 
If $(\Gam, \mu)$ has natural weights then $R(x,y) \le d(x,y)$, 
and if in addition $\Gam$ is a tree then $R(x,y)=d(x,y)$. 

The following is an easy consequence of (3.3).

\proclaim Lemma 3.1. 
For all $f\in \bR^{\Gamma}$ and $x,y\in \Gamma$,
$$|f(x)-f(y)|^2\le R(x,y)\sE(f,f). \eqno(3.4)$$
Further, for each $x,y\in \Gamma$, there exists $f$ so that the equality
holds in (3.4).

We recall  some basic properties of Green kernels.
Let $Y_t^B$ be the continuous time random walk on $(\Gamma, \mu)$ killed
outside $B:=B_R(x_0,r)$, and  $q_t^B(x,y)$ be the
transition density of $Y_t^B$. 
The Green kernel $g_B(x,y)$ of $Y_t^B$ is defined by
$g_B(x,y)=\int_{0}^{\infty}q_t^B(x,y)dt$. Then $g_B(\cdot,\cdot)$ 
has the reproducing property that 
$$ \sE(g_B(x,\cdot), f)=f(x) $$ 
for all $f\in H^2$ such that $f|_{B^c}=0$.

Using this and the fact that $e_{B,x}(y):=g_B(x,y)/g_B(x,x)$ is 
the equilibrium potential for $R(x,B^c)$, we have
$$ R(x,B^c)^{-1}=\sE(e_{B,x},e_{B,x})=g_B(x,x)^{-1}, $$
so that 
$$ R(x,B^c)=g_B(x,x)=\int_{0}^{\infty}q_t^B(x,x)dt\qquad\forall
x\in \Gamma, B\subset \Gamma.\eqno (3.5)$$

\bigbreak \noindent
{\bf 4. Heat kernel estimates on graphs and trees }

\ms
Recall that for $x\in \Gamma$ and $r\ge 0$, we 
denote $V(x,r)=\mu (B(x,r))$.

\proclaim Theorem 4.1. Let $(\Gam,\mu)$ be a weighted graph and suppose
that the edge weights satisfy $\mu_{xy}\ge 1$ for all
$x$ and $y$. Then 
$$ q_{2rV(x,r)}(x,x) \le \frac{2}{V(x,r)}, \quad x\in \Gamma, \,r>0.  $$ 

\med {\bf Remark.} This is similar to the bound in
Proposition 3.2 of [BCK], but has weaker hypotheses: in particular the
bound on $q_t(x,x)$ only uses the volumes of the balls $V(x,R)$.

\proof
Fix $x_0\in \Gamma$, write $B(r)=B(x_0,r)$ and  $V(r)=V(x_0,r)$. Set 
$f_t(y)=q_t(x_0,y)$ and 
$$ \psi(t) = ||f_t||_2^2= q_{2t}(x_0,x_0) = f_{2t}(x_0); $$
note that $\psi$ is decreasing. 
Let $r>0$; since 
$$ \sum_{y\in B(r)} f_t(y) \mu_y \le 1,$$ 
there exists $y=y(t,r) \in B(r)$ with $f_t(y)\le V(r)^{-1}$. 
Note that, since $\mu_{e}\ge 1$ for every edge $e$, it follows that 
$R(x,y)\le d(x,y)$ for all $x$, $y$. 
Then by (3.4)
$$ \eqalign{
\half f_t(x_0)^2 &\le  f_t(y)^2 + |f_t(x_0)-f_t(y)|^2   \cr
 &\le \frac{1}{V(r)^2} + R(x_0,y) \sE(f_t,f_t)
 \le \frac{1}{V(r)^2} + r \sE(f_t,f_t). }$$  
Hence
$$ \psi'(t) = -2 \sE(f_t,f_t) \le 
\frac{ 2 V(r)^{-2} - \psi(t/2)^2}{r}. \eqno(4.1) $$
Since $-\psi(s/2) \le -\psi(t)$ for $t\le s \le 2t$,
integrating (4.1) from $t$ to $2t$ we obtain
$$ \psi(2t)-\psi(t) \le  2t r^{-1} V(r)^{-2} - t r^{-1} \psi(t)^2. $$
So as $\psi(2t)>0$,
$$ t V(r)^2  \psi(t)^2 \le 2t + r  V(r)^2 \psi(t)\le (4t)\vee (2rV(r)^2 \psi(t)). $$
Hence
$$   \psi(t) \le \frac{2}{V(r)} \vee \frac{2r}{t}. $$
Taking $r$ such that $t=r V(r)$ completes the proof. \qed

\proclaim Corollary 4.2. Let $V(x,r)\ge r^2/A$, and $t=r^3$. Then
$$ q_{2t}(x,x) \le \frac{2 (A \vee 1)}{r^2}=
   \frac{2 (A \vee 1)}{ t^{2/3}}\, .  \eqno(4.2)$$

\proof Let $\lam=r^{-2} V(x,r)$, so that $\lam \ge A^{-1}$.
Let $t_0=r V(x,r)=\lam r^3$. 
If $\lam\le 1$ then $t_0\le t$ and so Theorem 4.1 gives 
$$ q_{2t}(x,x) \le q_{2t_0}(x,x) \le \frac{2}{V(x,r)}
 = \frac{2}{\lam r^2} =2 \lam^{-1} t^{-2/3}
  \le 2A  t^{-2/3}. $$ 

Now suppose that $\lam\ge 1$. Let $r'$ be such that
$t=r'V(x,r')$; as $rV(x,r)=\lam r^3=\lam t$, we have $r'\le r$.
So
$$  q_{2t}(x,x) = q_{2r'V(r')}(x,x) \le \frac{2}{V(x,r')}
 =  \frac{2r'}{t} \le  \frac{2r}{t}= 2 t^{-2/3}
  \le 2  (A \vee 1) t^{-2/3}. $$
\qed

\proclaim Lemma 4.3. Let $f_t(y)=q_t(x_0,y)$. Then
$$ \Big| \frac{f_t(y)}{f_t(x_0)}-1 \Big|^2 
 \le  \frac{d(x_0,y)}{t f_t(x_0)} .  \eqno(4.3)$$

\proof Let $e(t) =\sE(f_t,f_t)$. Then $e$ is decreasing, and
$$ |f_t(x_0)-f_t(y)|^2 \le d(x_0,y) e(t). $$
So as
$$ \psi(t)-\psi(t/2) =-2 \int_{t/2}^t e(s) ds, $$
we have
$$ 2 e(t)\cdot t/2 \le  2 \int_{t/2}^t e(s) ds
\le \psi(t/2). $$
So,
$$ |f_t(x_0)-f_t(y)|^2 \le \frac{d(x_0,y) f_t(x_0)}{t}, $$
and dividing by $f_t(x_0)^2$ completes the proof. \qed 

\ms Up to this point we have not needed to use the fact that
$\Gam$ is a tree, but the following lemma relies strongly
on this. From now on we take $\Gam$ to be a subgraph of $\bB$, 
and define $M(x,r)$, and the conditions $\lam$--good,
 $G_2(N,R)$ and $G_3(x,y,m,\kappa)$ as in Section 2.

\proclaim Lemma 4.4. Let $B=B(x_0,r)$, and $x\in B(x_0,r/8)$. 
Then
$$  \frac{r}{8 M(x_0,r)} \le g_B(x,x)=R(x,B^c)\le 9r/8. \eqno(4.4)$$

\proof Since $x$ is connected to $B(x,r)^c$ by a path of 
length $9r/8$, the upper bound is clear.  

For the lower bound let $m=M(x_0,r)$ and
$A=\{z_1, \dots, z_m\}$ be the set given in Definition 2.9:
note that $d(x,z_i)\ge r/8$ for each $i$.
Let $h_i$ be the function on $G$ such that
$h_i(z_i)=1, h_i(x)=0$ and $h_i$ is harmonic $G-\{x,z_i\}$.
Then $h_i(y)=\bP^y(T_{z_i}<T_x)$, and 
$$ \sE(h_i,h_i) = R(x,z_i)^{-1} = d(x,z_i)^{-1}\le \frac{8}{r}. $$
If $y\in B(x,r)^c$ then since any path from $y$ to $x$
passes through $A$, we have $h_i(y)=1$ for at least one $i$.
So if $h=\max_i h_i$ then $h(x)=0$ and $h= 1$ on $B(x,r)^c$.
So, using (3.1), 
$$ R(x,B^c)^{-1} \le \sE(h,h) \le m \max_i \sE(h_i,h_i)
  \le \frac{8M(x_0,r)}{r}, $$
proving the lower bound \qed 

\proclaim Lemma 4.5. Let $B=B(x_0,r)$, $M=M(x_0,r)$.
\nl (a) 
$$ E^z \tau_B \le 2r V(x_0,r), \qq z\in B(x_0,r). \eqno(4.5)$$
(b) 
$$  E^x  \tau_B \ge \frac{r V(x_0,r/(32M))}{32 M},
 \hbox{ for } x \in B(x_0, r/(32M)). \eqno(4.6)$$

\proof For any $z\in B$,
$$  E^z\tau_B = \sum_{y\in B} g_B(z,y)\mu_y. \eqno(4.7)$$
The upper bound follows easily from (4.7), since
$$  \sum_{y\in B} g_B(z,y)\mu_y \le \sum_{y\in B} g_B(z,z) \mu_y
 = R(z,B^c) V(x,r) \le 2 r V(x,r). $$

For the lower bound, let $x \in B(x_0,r/8)$, and 
set $p_B^x(y)=g_B(x,y)/g_B(x,x)$. Then
$\sE(p^x_B,p^x_B)= g_B(x,x)^{-1}$ and so 
$$ | 1-p^x_B(y)|^2 \le d(x,y) R(x,B^c)^{-1} \le d(x,y) (8M/r). $$
Let $B'= B(x_0,r/(32M))$. Then if $x, y \in B'$, 
$d(x,y) \le r/(16M)$ and so $p^x_B(y)\ge 1 - 2^{-1/2}\ge \fract14$.
So, using Lemma 4.4,
$$  E^x \tau_B \ge \sum_{y\in B'} g_B(x,x)p^x_B(y) 
  \ge \fract14 \mu(B') R(x,B^c) 
  \ge r \mu(B')/(32 M). $$
\qed

\proclaim Proposition 4.6. 
Let $r\ge 1$ and $x_0\in \Gamma$, and $B=B(x_0,r)$. 
Write $M=M(x_0,r)$, $V=V(x_0,r)$ and 
let $V_1=V_1(x_0,r)=V(x_0,r/(32 M(x_0,r)))$.
Then if $x \in B(x_0, r/(32M))$, 
$$  P^x(\tau_B \le t) \le \Big(1- \frac{V_1}{64MV}\Big) + \frac{t}{2rV}.  $$
and 
$$q_{2t}(x,x)\ge \frac{c_1 V_1(x_0,r)^2}{V(x_0,r)^3M(x_0,r)^2} \q
\hbox{ for } t\le  \frac{rV_1(x_0,r)}{64M(x_0,r)}. $$

\proof The proof is standard. By the Markov property, 
$$ \bE^{x}[\tau_B] \le t+\bE^x[1_{\{\tau_B> t\} } E^{Y_t}(\tau_B)] ,$$
for all $t>0$. Using this and Lemma 4.5,
$$ \frac{rV_1}{32 M} \le t + P^x(\tau_B>t) 2r V, $$
and rearranging this we have
$$ P^x( Y_t \in B) \ge  P^x(\tau_B>t) \ge \frac{ (rV_1/32 M)-t}{2r V}. 
\eqno(4.8) $$
This proves the first assertion.

By (4.8) if $t \le rV_1/(64M)$ then
$$ P^x( Y_t \in B) \ge \frac{c_2 V_1}{VM}. $$
By Chapman-Kolmogorov and Cauchy-Schwarz 
$$  P^x( Y_t \in B)^2 = ( \sum_{y \in B} q_t(x,y) \mu_y)^2
\le \mu(B) \sum_{y \in B} q_t(x,y)^2 \mu_y  \le q_{2t}(x,x) V. $$
So 
$$ q_{2t}(x,x) \ge V^{-1}  P^x( Y_t \in B)^2
 \ge  \frac{c_2^2 V_1^2}{ V^3M^2}. \eqno(4.9) $$
\qed

\proclaim Theorem 4.7. Suppose that $B=B(x_0,r)$ is $\lam$--good for 
$\lambda\ge 1$, and let  $I=I(\lam,r)=[r^3 \lam^{-6}, r^3 \lam^{-5}]$.
\nl (a) For $x \in B(x_0, r/\lam)$,
$$  c_0 \frac{r^3}{\lam^5} \le E^x \tau_B \le 2 \lam {r^3}. \eqno(4.10) $$
(b) For each $K \ge 0$ 
$$   q_{2t}(x_0,y) \le  (1+{\sqrt K}) t^{-2/3} \lam^{3}
\q \hbox{ for } t\in I, \q y\in B(x_0,Kt^{1/3}). \eqno(4.11)$$
(c) Let $x \in B(x_0, r/\lam)$. 
Then 
$$ q_{2t}(x,y) \ge c_1 t^{-2/3} \lam^{-17}, \q
 \hbox{ if } d(x,y) \le  c_2 \lam^{-19}r, \q    t \in I.  \eqno(4.12) $$

\proof (a) Let $B$, $V$, $V_1$, $M$ be as in the previous proof.
As $32M\le 64 M\le \lam$, $V_1\ge V(x,r/\lam) \ge r^2 \lam^{-4}$,
while $V\le \lam r^2$. Thus (4.10) is immediate from Lemma 4.5.

\noi (b) Let $t_1=(r/\lam^2)^3$. Then by Corollary 4.2 (taking $A=\lambda^2$), 
if $t\in I$,
$$ q_{2t}(x_0,x_0) \le  q_{2t_1}(x_0,x_0) 
 \le 2 \lam^2 t_1^{-2/3} \le 2 \lam^{8/3} t^{-2/3} 
\le \lam^3 t^{-2/3}. \eqno(4.13) $$
Now, 
for $t\in I$ and $y\in B(x_0,Kt^{1/3})$, we have, using Lemma 4.3
and (4.13),
$$\eqalign{ 
q_{2t}(x_0,y)&\le  q_{2t}(x_0,x_0)+| q_{2t}(x_0,y)-q_{2t}(x_0,x_0)|\cr
&\le q_{2t}(x_0,x_0)+\sqrt{\frac K{2t^{2/3}}q_{2t}(x_0,x_0)}\le 
(1+{\sqrt K}) t^{-2/3} \lam^{3},\cr }$$
proving (4.11).

\noi (c) Let $x \in B(x_0, r/\lam) \subset B(x_0, r/(32M))$. 
Then $rV_1/(64M) \ge r^3 \lam^{-5}$, so for $t\in I$
by Proposition 4.6,
$$ q_{2t}(x,x) \ge c_2 V_1^2 /(V^{3}M^{2}) \ge c_2 r^{-2} \lam^{-13}
 \ge c_2 t^{-2/3} \lam^{-17}, $$
where $c_2=c_{4.6.1}$.
Hence, by Lemma 4.3, if $d(x,y) \le c_2 \lam^{-19}r$,
$$ \Big| \frac{q_{2t}(x,y)}{q_{2t}(x,x)}-1 \Big|^2 
        \le  \frac{d(x,y)}{2t q_{2t}(x,x)}
\le \frac{d(x,y) r^2\lam^{13}}{2c_2 t}
\le \frac{d(x,y)\lam^{19}}{2c_2 r}\le \half,$$
from which (4.12) follows. \qed

\proclaim Corollary 4.8. Let $\lam\ge 64$, 
and $B(x,r)$ and $B(x, \lam^{-5} r)$ be $\lam$--good.
Then
$$ E^x d(x,Y_t) \ge c_1 \lam^{-12} t^{1/3}, 
 \q \hbox { for } 
\frac{r^3}{\lam^6} \le t \le \frac{r^3}{\lam^5}. $$

\proof Let $I=[r^3 \lam^{-6}, r^3 \lam^{-5}]$ and 
$B'=B(x,r \lam^{-5})$.
Let $t\in I$, and $y\in B'$. Then since $r\le \lam^2 t^{1/3}$,
$d(x_0,y) \le \lam^{-5} r \le \lam^{-3} t^{1/3}$, so by
(4.11) (with $K=1)$ we have $q_{2t}(x_0,y) \le 2 t^{-2/3} \lam^3$.
Hence since $B'$ is $\lam$--good,
$$  P^x( Y_{2t} \in B')= \sum_{y\in B'} q_{2t}(x_0,y) \mu_y
 \le \mu(B')  2 t^{-2/3} \lam^3 \le 2 \lam^{-2} \le \half. $$
Thus
$$ E^x d(x,Y_{2t}) \ge  \lam^{-5} r P^x(Y_{2t} \not\in B')
   =  \lam^{-5} r (1- P^x(Y_{2t} \in B'))\ge \half r \lam^{-5}. $$
\qed

\proclaim Lemma 4.9. Suppose $x$ satisfies $G_2(N,R)$. Then
$$ P^x( \tau_{B(x,NR)} \le t ) \le e^{-c_1 N} \q 
\hbox{ provided  $N \ge c_2 t/R^3$.} $$

\proof We use the argument of [BB1]. Let 
$$ A = \{ y \in G: B(y,R/2) \hbox{ is $\lam_1$--good} \} . $$
Define stopping times $(T_i)$, $(S_i)$ by taking 
$T_0=\min\{ t: Y_t \in A\}$, and
$$ \eqalign {
 S_{n} &=\min \{ t \ge T_{n-1}: Y_t \not\in B(Y_{T_{n-1}},R/2) \}, \cr
 T_{n} &= \min \{ t \ge S_n: Y_t \in A \} .} $$
Since $x$ satisfies $G_2(N,R)$ we have $T_{N/8} \le  \tau_{B(x,NR)}$ 
$P^x$ -a.s.
Let $\xi_i=S_{i+1}-T_{i}$, $i\ge 1$. Then by Proposition 4.6
there exists $p=p(\lam_1)<1$ and $c_3=c_3(\lambda)>0$ such that 
$$ P^x \big(\xi_i \le s| \sigma(Y_u, 0\le u \le T_i )\big) 
 \le p + c_3 R^{-3} s .  \eqno(4.14) $$
Lemma 1.1 of [BB1] (see also Lemma 3.14 of [B1]) gives that, writing
$a=c_3/R^3$, (4.14) implies that
$$ \log P^x( \sum_{i=1}^{N/8} \xi_i \le t ) \le
  -\fract18 N \log(1/p) + 2 \Big( \frac{aNt}{8p} \Big)^{1/2}. $$
Substituting for $a$ we deduce that
$$  \log P^x( \tau_{B(x,NR)} \le t ) 
 \le  -N\Big( 2c_4 -c_5( t/(R^3N))^{1/2}  \Big) \le {-c_4N}, $$
provided $N \ge (c_5/c_4)^2\cdot (t/R^3)$. \qed

\proclaim Theorem 4.10. 
Let $x,y \in \sG$, $t>0$ be such that $N:=[\sqrt{d(x,y)^3/t}]\ge 8$  
and suppose the event 
$F_*(x,y,d(x,y)N^{-1},\frac 18 N; d(x,y)^3t^{-2/3},N)$ holds. Then 
$$ q_t(x,y) \le c_1 t^{-2/3} \exp( -c_2 N).  \eqno(4.15)$$

\proof Define $T_{z_0}=\inf \{t:Y_t={z_0}\}$ and $R=d(x,y)/N$, 
where $z_0$ is a middle point in $\gamma (x,y)$. 
Let $G_x$ be the set of points $w$ in $\sG$ such that $\gamma(x,w)$ does not
contain $z_0$, and let $G_y=\sG-G_x$. 
Then, we have
$$\eqalignno{
  q_{t}(x,y)\mu_x \mu_y &= \mu_x P^x(Y_t=y) \cr
  &= \mu_x P^x(Y_{t/2} \in G_y, Y_t=y)+ \mu_x P^x(Y_{t/2} \in G_x, Y_t=y)\cr
&=  \mu_x P^x(Y_{t/2} \in G_y, Y_t=y)+ \mu_y P^y(Y_{t/2} \in G_x, Y_t=x),
 &(4.16) \cr}$$
where in the last line we used the $\mu$--symmetry of $Y$.
The two terms in (4.16) are bounded in the same way. For the first,
$$\eqalignno{
 P^x(Y_{t/2} \in G_y, Y_t=y) &\le  P^x( T_{z_0}\le t/2, Y_t=y) \cr
 &= E^x \big( 1_{( T_{z_0}\le t/2)} P^{z_0}(Y_{t-T_{z_0}}=y) \big) \cr
 &\le P^x(T_{z_0}\le t/2)  \sup_{t/2\le s\le t} q_s(z_0,y)\mu_y.\cr
 &\le \mu_y \sqrt{q_{t/2}(y,y)q_{t/2}(z_0,z_0)}P^x(T_{z_0}\le t/2) \cr
 &\le \mu_y N^3t^{-2/3} P^x(T_{z_0}\le t/2), }$$
where we used (4.11) with $\lambda=N, r=N^2t^{1/3}$ 
in the last inequality. 
Now, $t/R^3\sim (d(x,y)^3/t)^{1/2}\sim N$, so $N\ge ct/R^3$. 
Thus, by Lemma 4.9 we have 
$$ P^x(T_{z_0}\le t/2)\le e^{-cN}~~\hbox{ and }~~P^y(T_{z_0}\le t/2)\le e^{-cN}.$$ 
Combining these facts 
$$ q_{t}(x,y)\le c' N^3t^{-2/3} e^{-cN} \le c t^{-2/3} e^{-c''N},$$
which completes the proof.\qed

\proclaim Theorem 4.11. Let $x,y \in \sG$, $m\ge 1$, $\kappa \ge 1$  and 
suppose $G_3(x,y,m,\kappa)$ holds. Then if $T= d(x,y)^3\kappa /m^2$
$$ q_{2T}(x,y) \ge c_1 T^{-2/3} e^{ -c_2 (\kappa+c_3) m}.  \eqno(4.17)$$

\proof Let $r=d(x,y)/m$, and $(z_i)$, $(\Th_i)$ be the points and integers
given by the condition  $G_3(x,y,m,\kappa)$ in Definition 2.15.
Let $B_i=B(z_i, \Th_i^{20}r)$, and $B'_i=B(z_i, r)$. 
Applying (4.12) to $B_i$ we deduce that if 
$d(y,y')\le c_{4.7.2} \th^{-19}(\Th_i^{20} r)$, and
$$ \Th_i^{54} r^3 \le t_i \le \Th_i^{55} r^3,  \eqno(4.18) $$
then
$$ q_{2t_i}(y,y') \ge c_4 t_i^{-2/3} \Th_i^{-17}. \eqno(4.19)$$
If $y_i \in B_i'$ then by the choice of $\Th_i$
$$ d(y_{i-1},y_i) \le 3 r \le  c_{4.7.2} \Th_i^{-19} (\Th_i^{20} r), $$
and so the bound in (4.19) holds for $q_{2t_i}(y,y')$.
Therefore for $y_{i-1} \in B_{i-1}'$ and $t_i$ satisfying (4.18),
$$ \int_{B'_i} q_{2t_i}(y_{i-1},y_i) \mu(dy_i)
 \ge c_4 t_i^{-2/3} \Th_i^{-17} \mu(B'_i)   \ge c_4 \Th_i^{-c_5}; $$
we used here the fact that $\mu(B'_i)\ge \Th_i^{-2} r^2$.
So if $t_i$ satisfy (4.18), and $s=\sum t_i$ then since
$\sum \log \Th_i \le \sum \Th_i^{54} \le m \kappa$, 
$$\eqalign{
&  q_{2s}(x,y) \ge
 \int_{B'_1} \dots \int_{B'_{m-1}} q_{2t_1}(x,y_1) q_{2t_1}(y_1, y_2) \dots
  q_{2t_{m}}(y_{m-1}, y) \mu(dy_1)\dots  \mu(dy_{m-1}) \cr
 &\ge  (c t_{m}^{-2/3} \Th_m^{-17}) c_4^{m-1} \Pi_{i=1}^{m-1} \Th_i^{-c_5} 
\ge s^{-2/3} \exp (-c_6m -c_5 \sum \log \Th_i )\cr
&\ge s^{-2/3} e^{ -(c_5\kappa+c_6)m} . \cr }$$

As  $G_3(x,y,m,\kappa)$ holds we have $r^3 \sum \Th_i^{54}\le m \kappa r^3 =T$. 
If $T\le  r^3 \sum \Th_i^{55}$ we can choose $(t_i)$ satisfying (4.18) so 
that $s=T$. If not, let $s'=T-s$, so that $s'\le  m \kappa r^3$.
Fix a $j$ such that $\Th_j$ is minimal and in
the chaining argument above add $m'$ extra steps 
(of time length $t'$ satisfying (4.18) for $i=j$) between $B'_{j-1}$ and $B'_j$.
Since $c_7^{54} \le \Th_j^{54} \le \kappa$, we have $c_8 r^3 \le t'\le \kappa r^3$.
Then choose $m', t'$ so that $m' t' + s=T$; we have $m'\le c m$. 
Each extra step gives a factor of $c_4 \Th_j^{-c_5}$ in the lower
bound in the chaining argument, so the total contribution multiplies the
lower bound by a number greater than $e^{-c (\kappa+c') m}$. Thus (4.17) holds. \qed

\bigb {\bf 5. Random walk on the conditioned critical GW-branching precess  } 

\ms
In this section, we state and prove our main results on the random
walk on the IIC. As in Section 2 we write $\sG$ for the IIC on 
$\bB$, and $\tbP$ for its law. 
Let $Y=\{Y_t\}_{t\ge 0}$ be the simple random walk on $\sG(\om)$ defined in
Section 3; we write $E_\om^x$ for its law of $Y$ started at $x$.  Let
$q^{\omega}_t(x,y)$ be the transition density of $Y$. 

\med {\it Proof of Theorem 1.2.}
Fix $x \in \bB$, and let $c_3=c_{2.12.2}$.
Let $a=2/c_3$ and $\lam_n=e+a \log n$,
and $r_n$ satisfy $r_n^3 \lam_n^{-6} = e^n$.
Let $F_n$ be the event that $B(x, r_n)$ is $\lam_n$--good. Then
by Corollary 2.12
$$ \tbP (F_n^c) \le c e^{-c_3 a \log n}= c' n^{-2}, $$
so by Borel-Cantelli $F_n^c$ occurs for only finitely many $n$, 
$\tbP$--a.s.
Let $N$ be the largest $m$ such that $F_m^c$ occurs; then 
$$ \tbP(N > m) \le \sum_{m+1}^\infty   \tbP (F_n^c) \le c m^{-1}. $$ 
Set $S(x)=e^N$.
For $n\ge (\log S(x))+1$ 
we have, by (4.11) and (4.12),
$$ c' t^{-2/3} \lam_n^{-17} \le q_{2t}(x,x) \le c'' t^{-2/3} \lam_n^3 \eqno(5.1)$$
for $e^n \le t\le \lam_n e^n$. 
Let $n(t)$ be the unique integer such that $\log t \in [n(t)-1, n(t))$.
Hence, if $t\ge S(x)$, $n(t)> N$ and so (5.1) holds for $n=n(t)$.
Since 
$$ \lam_{n(t)} = e+ a \log n(t) \sim a \log\log t, $$
we obtain (5.1). \qed

\sms While the powers of the terms in $\log\log t$ 
given in Theorem 1.2 are not the best 
possible, we do have oscillations in $t^{-2/3}q^\om_t(.,.)$ of that order.

\proclaim Lemma 5.1.
$$ \liminf_{t \to \infty} 
   (\log\log t)^{1/6} t^{2/3} q^\om_{2t}(0,0)  \le 2, 
\q P_\om^0 -a.s. \eqno(5.2) $$

\proof Define $a_n$ by $V(0,2^n)= a_n 2^{2n}$, and
let $t_n = 2^n V(0,2^n) =  a_n 2^{3n}$. Then by Theorem 4.1,
$$ q^\om_{2t_n}(0,0) \le \frac{2}{V(0,2^n)}=  \frac{2t_n^{-2/3}}{a_n^{1/3}}.  $$
By Proposition 2.8(a), $a_n> (\log n)^{1/2}$ for infinitely many $n$, a.s.,
giving (5.2). \qed

\bigb {\it Proof of Theorem 1.3}.
(a) The lower bound in (1.4) is an immediate consequence of 
Corollaries 2.12 and 4.8.
For the upper bound, let $Z_t= \sup_{0\le s\le t} d(x,Y_s)$,
$R=t^{1/3}$ and $T_M=\tau_{B(x,MR)}$.
Let $K_t(x)(\om)$ be the largest $n$ such that $x$ does not 
satisfy $G_2(n,R)$. Then by Proposition 2.14
$$ \tbP_x ( K_t(x) \ge k) 
\le \sum_{l=k}^\infty\tbP_x (x \hbox{ does not satisfy } G_2(l,R))
\le  c' e^{-ck}. \eqno(5.3) $$
Then $\{Z_t \ge nR \}\subset \{T_n\le t\}$, and so by Lemma 4.9, 
$$  \eqalignno{
 E^x_\om Z_t &\le R \sum_{n=0}^\infty  P^x_\om(T_n \le t) \cr
 &\le  R\Big( 1+K_t(x) + 
 \sum_{n=K_t(x)+1}^\infty  P^x_\om(T_n \le t)\Big) \cr
&\le R\Big(1+ K_t(x) +   \sum_{n=K_t(x)+1}^\infty  e^{-cn} \Big) 
\le R(c+  K_t(x)). &(5.4) \cr }$$
Since $ \tbE_x K_t(x) \le c'$ this completes the proof. 

\noi (b) Let $m(t) =\lfloor t \rfloor$; 
Since 
$$ |E^x_\om d(x,Y_t)- E^x_\om d(x,Y_{m(t)})|
 \le E^x_\om d(Y_{m(t)},Y_t) \le c, $$
it is enough to prove (1.5) for integer $t$.
Using (5.3) and Borel-Cantelli there exists $c'$ such that
$$ \tbP_x( K_n(x) > c' \log n ~~\hbox {i.o.}) =0. $$
and so by (5.4) 
$$ E^x_\om d(x,Y_n) \le c'' n^{1/3} \log n $$
for all sufficiently large $n$.
The lower bound in (1.5) follows from Corollary 4.8 by the same
argument as in Theorem 1.2.
\qed

\bigb {\it Proof of Theorem 1.4}.
We begin with the on-diagonal case $x=y$.
Let $\lam_n=n$ and $r_n$ be defined by $2r_n^3/\lam_n^6=t$.
Let $F_n=\{ B(x,r_n)$ is  $\lam_n$--good $\}$, 
and  $N(\om)=\min\{ n: \om \in F_n\}$. 
By Corollary 2.12 
$\tbP_x(N> n)\le \tbP_x(F_n^c) \le e^{-cn}$.
On $F_n$ we have, by (4.11), 
$ q^\om_t(x,x) \le c t^{-2/3} n^3, $
so 
$$  \tbE_x [q^{\omega}_t(x,x)]\le  c t^{-2/3}  \tbE_x N^3
 \le c'  t^{-2/3}, \eqno(5.5)$$ 
proving the on-diagonal upper bound in (1.6).

For the on-diagonal lower bound choose $m_0$ such that
$\tbP_x(F_{m_0})\ge \half$ and then on $F_{m_0}$, by the lower
bound in (4.12),
$$  q^\om_t(x,x) \ge  c t^{-2/3} m_0^{-17}. $$

\ms For the off-diagonal bounds,
when $d(x,y)\le 64t^{1/3}$, (1.6) can be proved similarly to 
(5.5) using Theorem 4.7(b). So we will assume $d(x,y)> 64t^{1/3}$.
Now, let $N:=[\sqrt{d(x,y)^3/t}]\ge 8$ and define 
$F_0=F_*(x,y,d(x,y)N^{-1},\frac 18 N; d(x,y)^3t^{-2/3},N)$.
Let $\lam_0=N$ and define $\lam_n=N+n-1$ for $n\ge 1$.
For each $n\ge 1$, set $r_n=t^{1/3}\lam_n^2$ and let
$F_n=\{B(x,r_n)$ is $\lam_n$-good $\}$. Then, $\tbP_{x,b}(F_n^c)\le e^{-c\lam_n}$. 
We now apply Theorem 4.7 (b) with
$K=\lam_n^2$ and obtain the following. 
(Note that we can apply the theorem because $d(x,y)/t^{1/3}\le cN^{2/3}\le c\lam_n^2$.) 
$$q_{2t}(x,y)\le c(1+\sqrt{\lam_n^2})t^{-2/3}\lam_n^3\le c't^{-2/3}\lam_n^4.\eqno(5.6)$$
Let $M(\omega)=\min\{n\ge 0: \omega\in F_n\}$. Then, 
$\tbP_x(M=0)=\tbP_x(F_0)\ge 1-4e^{-N}$ and 
$\tbP_x(M>n)\le \tbP_x(F_n^c)\le ce^{-c'\lam_n}$.
Thus, using Theorem 4.10 and (5.6), we obtain
$$ \eqalign{\tbE_{x,y} [q^{\omega}_t(x,y)]=&\tbE_{x,y} [q^{\omega}_t(x,y):M=0]
+\tbE_{x,y} [q^{\omega}_t(x,y):M>0]\cr
\le & ct^{-2/3}\exp(-c'N)+c''t^{-2/3}\tbE[\lam_M^4:M>0].}$$
Since $\tbE[\lam_M^4:M>0]\le c\sum_{k=1}^\infty(N+k-1)^4e^{-c'(N+k-1)}\le ce^{-c''N}$,
we obtain (1.6).

We next prove (b). Choose $\kappa= 2 c_{2.16.1}$, so that 
$\bP_{x,y}( G_3(x,y,m,\kappa) \hbox{ holds })\ge \half$. 
Now choose $m=(R^3\kappa/t)^{1/2}$; by Theorem 4.11, for $\om$ such that 
$G_3(x,y,m,\kappa)$ holds, 
$$ q_{2t}^\om(x,y) \ge    c t^{-2/3} \exp (-c'(\kappa+c'') m). $$
Taking expectations gives (1.7). \qed

\ms Let 
$$ \wt Z^{(n)}_t = n^{-1/3} d(0,Y_{nt}), \q t\ge 0. $$
By Theorem 1.3(a) the process $\wt Z^{(n)}$ is tight
with respect to the annealed law given by the
semi-direct product $\bP^*= \bP \times P^0_{\omega}$. 
(See Theorem 1.1 for the analogous result for the 
discrete time simple random walk.)

\med {\it Proof of Theorem 1.5.}
Let $U_n= \sup_{0\le s \le 1}  Z^{(n)}_s$. Then, by (4.5), 
$$ \eqalign{
 P^0_{\omega}( U_n \le \lam) &= 
  P^0_{\omega}( \sup_{t\le n} d(0,Y_s) \le \lam n^{1/3}) \cr
 &= P^0_{\omega}( \tau_{B(0,\lam n^{1/3})}  \ge n ) 
\le \frac{ 2 \lam n^{1/3} V(0,\lam n^{1/3})}{n}  .} $$
So by Proposition 2.8(b), we have, for any $\lam>0$, that  
$\liminf_{n\to\infty}  P^0_{\omega}( U_n \le \lam) =0$, which shows that
the r.v. $U_n$ (and hence the processes  $Z^{(n)}$) are not
tight. \qed

\med {\bf Remark.}
 This result illustrates the difference in the type of results
that can arise between the quenched and annealed cases. 
For the case of supercritical bond percolation in $\bZ^d$,
while an invariance principle was proved
in the annealed case  in [DFGW] in 1989, the quenched case
(for $d\ge 4$) was only proved recently in [SS]. 

\medskip
\noindent {\bf Acknowledgment. }
The authors thank 
Antal J\'arai, Harry Kesten and Gordon Slade for valuable comments. 
%Harry Kesten for valuable comments. 

\bigskip \centerline {\bf References}
\def\nr{\noindent}

%\nr [Ai] M. Aizenman. On the number of incipient spanning clusters.
%{\it Nuclear Phys. B} {\bf 485} (1997), 551--582.

\nr [AO] S. Alexander, R. Orbach. Density of states on
fractals: \lq\lq fractons \rq\rq. {\it J. Physique (Paris) Lett.} 
{\bf 43}, L625--L631 (1982).

\noindent
[B1] M.T. Barlow.
Diffusions on fractals. 
Lectures in Probability Theory and Statistics:
Ecole d'\'et\'e de probabilit\'es de Saint-Flour XXV, Springer, 
New York, 1998.

\noindent [B2] M.T. Barlow.
Random walks on supercritical percolation clusters.
{\sl Ann. Probab.} {\bf 32} (2004), 3024-3084. 

\nr [BB1] M. T. Barlow, R. F. Bass.  The construction of Brownian motion 
on the Sierpinski carpet.  {\it Ann. Inst. Henri Poincar\'e \bf 25} (1989),
225--257.

\nr [BB2] M.T. Barlow, R.F. Bass. Random walks on graphical Sierpinski 
carpets. {\it Random walks and discrete potential theory (Cortona, 1997)}, 
26--55, Sympos. Math., XXXIX, Cambridge Univ. Press, Cambridge, 1999.

\noindent
[BCK] M.T. Barlow, T. Coulhon, T. Kumagai.
Characterization of sub-Gaussian heat kernel estimates on strongly 
recurrent graphs. To appear {\it Comm. Pure Appl. Math.}

\nr [BCKS] 
C. Borgs, J.T. Chayes, H. Kesten, J. Spencer. 
The birth of the infinite cluster: finite-size scaling in percolation. 
{\it Comm. Math. Phys.}  {\bf 224} (2001),  no. 1, 153--204.

\noindent
[DFGW] A. De Masi, P.A. Ferrari, S. Goldstein, W.D. Wick. 
An invariance principle for reversible Markov processes. 
Applications to random motions in random environments. 
{\it J. Stat. Phys.} {\bf 55} (1989), 787--855. 

\noindent [F] D. Freedman. On tail probabilities for martingales.
{\it Ann. Probab.} {\bf 3} (1975), 100-118.

\nr [GT1] A. Grigor'yan, A. Telcs. Sub-Gaussian estimates of heat kernels
on infinite graphs. {\it Duke Math. J.} {\bf 109} (2001), 452-510.

\nr [GT2] A. Grigor'yan, A. Telcs. Harnack inequalities and sub-Gaussian 
estimates for random walks.  {\it Math. Annalen} {\bf 324} (2002), 521--556. 

\noi [Gm] G.R. Grimmett. {\it Percolation.} (2nd edition). Springer, 1999.

\noi [HK] B.M. Hambly, T. Kumagai. 
Heat kernel estimates for symmetric random walks on a class of fractal graphs
and stability under rough isometries. 
{\it Fractal geometry and applications: A Jubilee of B. Mandelbrot (San Diego, CA, 2002)},  
233--260, Proc. Sympos. Pure Math., 72, Part 2, Amer. Math. Soc., Providence, RI, 2004.

\noi [HS] T. Hara, G. Slade.
Mean-field critical behaviour for percolation in high dimensions.  
{\it Comm. Math. Phys.} {\bf 128} (1990),  no. 2, 333--391. 

\noindent [Har] T.E. Harris. 
The theory of branching processes. 
Dover Publications, Inc., New York, 2002. 
(Originally; Springer-Verlag, Berlin, 1963).

\noindent [vH] R. van der Hofstad.
Infinite canonical super-Brownian motion and scaling limits.
Preprint 2004.

\nr [HJ] R. van der Hofstad, A.A. J\'arai. 
The incipient infinite cluster for high-dimensional unoriented percolation.  
{\it J. Stat. Phys.} {\bf 114} (2004), 625-663.

\nr [HHS] R. van der Hofstad, F. den Hollander, G. Slade.
Construction of the incipient infinite cluster for 
spread-out oriented percolation above $4+1$ dimensionals.
{\sl Comm. Math. Phys.} {\bf 231} (2002), 435-461.

\nr  [HS] R. van der Hofstad, G. Slade.
Convergence of critical oriented percolation to super-Brownian motion 
above $4+1$ dimensions.  
%{\it Ann. Inst. H. Poincaré 
{\it Ann. Inst. Henri Poincar\'e 
Probab. Statist.} {\bf 39} (2003), no. 3, 413--485. 

\nr [Ja1] A.A. J\'arai.
Incipient infinite percolation clusters in 2D.  
{\it Ann. Probab.} {\bf 31}  (2003),  no. 1, 444--485.

\nr [Ja2] A.A. J\'arai.
Invasion percolation and the incipient infinite cluster in 2D.  
{\it Comm. Math. Phys.} {\bf 236}  (2003),  no. 2, 311--334. 

\nr [Jo] O.D. Jones. Transition probabilities for the simple random
walk on the Sierpinski graph. {\it Stoch. Proc. Appl. \bf 61} (1996),  45-69.

\nr [Ke1]  H. Kesten. The incipient infinite cluster in two-dimensional
percolation. {\it Probab. Theory Related Fields} {\bf 73} (1986), 369--394. 

\noindent
[Ke2] H. Kesten.
Subdiffusive behavior of random walk on a random cluster.
{\it Ann. Inst. Henri Poincar\'e} {\bf 22} (1986), 425--487.

\nr [Ke3] H. Kesten. Subadditive behavior of random walk on a
random cluster. Unpublished notes. 

\noindent
[Kig] J. Kigami. 
{\sl Analysis on Fractals.}
Cambridge Univ. Press, Cambridge, 2001.

\noindent  
[SS] V. Sidoravicius, A.-S. Sznitman. Quenched invariance principles
for walks on clusters of percolation or among random conductances.
{\it Probab. Theory Related Fields} {\bf 124} (2004), 219--244.

\vskip 0.2 truein

\noi{\ver, \q \vdate}
\vskip 0.2 truein

\noi MTB: Department of Mathematics, University of British Columbia, Vancouver
V6T 1Z2, Canada
\ms
\noi TK: Research Institute for Mathematical 
Sciences, Kyoto University, Kyoto 606-8502, Japan

\end